\documentclass[preprint,12pt]{elsarticle}

\usepackage[utf8]{inputenc}    
\usepackage[T1]{fontenc}       
\usepackage{mathpazo}          
\usepackage{beramono}          

\usepackage{ifxetex}
\ifxetex
    \usepackage{fontspec}
    \newfontfamily\zhawtitlefont{Helvetica Rounded Bold}
\else
    \newcommand{\zhawtitlefont}{\scshape}
\fi

\usepackage{xcolor}
\definecolor{zhawblue}{rgb}{0.00, 0.39, 0.65}
\definecolor{codegreen}{rgb}{0,0.6,0}
\definecolor{codegray}{rgb}{0.5,0.5,0.5}
\definecolor{codepurple}{rgb}{0.58,0,0.82}
\definecolor{codebackground}{rgb}{0.93,0.94,0.95}
\definecolor{darkgreen}{rgb}{0.0, 0.5, 0.0}
\definecolor{darkred}{rgb}{0.65, 0.0, 0.0}
\definecolor{darkblue}{rgb}{0.0, 0.0, 0.65}
\definecolor{darkyellow}{rgb}{0.75, 0.6, 0.0}
\definecolor{custompurple}{RGB}{148, 103, 189}
\definecolor{darkorange}{RGB}{180, 120, 0}

\usepackage{amsmath, amssymb}
\usepackage{algorithm}
\usepackage{algpseudocode}
\usepackage{amsthm} 

\allowdisplaybreaks

\usepackage{booktabs}
\usepackage{array}
\usepackage{tabularx}
\usepackage{multirow}
\usepackage{makecell}
\usepackage{adjustbox}

\usepackage{graphicx}
\usepackage{caption}
\usepackage{subcaption}
\usepackage{tikz}
\usepackage{pgfplots}
\pgfplotsset{compat=1.18}
\usetikzlibrary{arrows.meta, decorations.pathreplacing, patterns, pgfplots.fillbetween}
\usepackage{pifont}

\usepackage{listings}
\lstdefinestyle{mystyle}{
    backgroundcolor=\color{codebackground},   
    commentstyle=\color{codegreen},
    keywordstyle=\color{magenta},
    numberstyle=\tiny\color{codegray},
    stringstyle=\color{codepurple},
    basicstyle=\ttfamily\footnotesize,
    breaklines=true,
    keepspaces=true,
    numbers=left,
    numbersep=5pt,
    showspaces=false,
    showstringspaces=false,
    showtabs=false,
    tabsize=4
}
\usepackage[colorlinks=true, allcolors=blue]{hyperref}
\usepackage{url} 

\usepackage{cleveref}
\usepackage{pdfpages}          
\usepackage{filecontents}
\usepackage{todonotes}
\usepackage{pdflscape}

\journal{a journal}

\begin{document}

\begin{frontmatter}



\title{AI2STOW: End-to-End Deep Reinforcement Learning to Construct Master Stowage Plans under Demand Uncertainty}

\author[inst1]{Jaike van Twiller\corref{cor1}}\ead{jaiv@itu.dk}
\cortext[cor1]{}
\author[inst1]{Djordje Grbic}
\author[inst1]{Rune Møller Jensen}

\affiliation[inst1]{organization={IT University of Copenhagen},
            addressline={Rued Langgaards Vej 7}, 
            city={Copenhagen},
            postcode={2300}, 
            country={Denmark}}

\begin{abstract}
The worldwide economy and environmental sustainability depend on efficient and reliable supply chains, in which container shipping plays a crucial role as an environmentally friendly mode of transport. Liner shipping companies seek to improve operational efficiency by solving the stowage planning problem. Due to many complex combinatorial aspects, stowage planning is challenging and often decomposed into two NP-hard subproblems: master and slot planning. This article proposes AI2STOW, an end-to-end deep reinforcement learning model with feasibility projection and an action mask to create master plans under demand uncertainty with global objectives and constraints, including paired block stowage patterms. Our experimental results demonstrate that AI2STOW outperforms baseline methods from reinforcement learning and stochastic programming in objective performance and computational efficiency, based on simulated instances reflecting the scale of realistic vessels and operational planning horizons.
\end{abstract}


\begin{highlights}
\item Scalable algorithms are crucial for solving complex stowage planning problems.
\item AI2STOW, an end-to-end deep reinforcement learning model, constructs master plans under demand uncertainty, including paired block stowage patterns.
\item Evaluation on a realistic large-scale vessel and operational planning horizons.
\item Efficient generation of adaptive, feasible solutions outperforming baseline methods from reinforcement learning and stochastic prgramming.
\item Strong empirical evidence supporting deep reinforcement learning as a scalable solution method.
\end{highlights}


\end{frontmatter}

\section{Introduction} \label{sec_drl_scale:intro}
Product availability and efficient deliveries depend on the smooth operation of complex supply chains to meet the demands of our dynamic and global societies. During the last century, maritime transport has emerged as the cornerstone of global trade and modern consumerism. About 45\% of annual transported goods, valued at \$8.1 trillion in global trade \cite{united_nations_conference_on_trade_and_development_review_2021}, are transported by container vessels. This substantial economic impact is accompanied by a significant environmental consequence, contributing to more than 400 million tonnes of CO2 emissions yearly \cite{lloyds_list_shipping_2022}.
To put this in perspective, a typical passenger vehicle emits about 4.6 tonnes of CO2 per year. However, it is worth noting that container vessels emit significantly less CO2 per cargo tonne-kilometer than other modes of transportation \cite{jensen_container_2018}.

The maritime transportation business is highly competitive; therefore, liner shipping companies offer global shipments at low prices. Consequently, profit margins are often slim and improving operational efficiency is of paramount importance. As a result, liner shipping companies create stowage plans to allocate containers to vessel capacity \cite{jensen_container_2018}. Stowage planning is acknowledged as a complex combinatorial optimization (CO) problem \cite{botter_stowage_1992, larsen_heuristic_2021, tierney_complexity_2014}, characterized by numerous interdependent objectives and constraints, some of which are NP-hard. This includes, but is not limited to, vessel capacity limitations, seaworthiness regulations, revenue maximization, and operational cost minimization, all of which must be addressed under conditions of cargo demand uncertainty.

Despite its substantial economic and environmental significance, stowage planning remains underexplored, particularly when compared to more established domains such as vehicle routing~\cite{jensen_container_2018, van_twiller_literature_2024}. The scarcity of publications suggests a relatively immature field, with opportunities for advancement in problem modeling, standardized benchmark datasets, and algorithmic evaluation~\cite{van_twiller_literature_2024}. Furthermore, many combinatorial aspects essential to solving representative real-world instances are often overlooked in existing research. The same holds for the problem's inherent stochasticity, which remains insufficiently addressed in the current literature \cite{van_twiller_literature_2024}.

In light of these limitations, prior research has often relied on simplifications to make the problem more tractable. Specifically, it can be observed that most contributions are unable to solve representative versions of the full stowage problem (e.g.,~\cite{botter_stowage_1992, larsen_heuristic_2021}). A common strategy to address this complexity is hierarchical decomposition into two sequential subproblems: the master planning problem (MPP) and the slot planning problem (SPP) (e.g.,~\cite{pacino_fast_2011, wilson_container_2000}). Despite various attempts, the search for scalable algorithms capable of solving these decomposed yet representative problems remains an open challenge.

In recent years, machine learning (ML), and in particular deep reinforcement learning (DRL), has shown significant potential to complement or enhance traditional CO techniques~\cite{bengio_machine_2021, mazyavkina_reinforcement_2021}. A wide range of CO problems has been effectively addressed using ML-based methods, which often excel in scalability, robustness to uncertainty, and adaptability to dynamic conditions \cite{kool_attention_2019,van_twiller_navigating_2025,sadana_survey_2025}. These ML4CO approaches enable the derivation of flexible, data-driven heuristics that are often infeasible to design manually. Nevertheless, the application of ML to stowage planning remains relatively limited. Existing approaches frequently fall short in capturing the full complexity of real-world stowage scenarios \cite{van_twiller_literature_2024}. This highlights the importance of a more in-depth investigation into ML-driven approaches specifically designed to address the unique challenges of stowage planning.

This article builds upon our previous work~\cite{van_twiller_navigating_2025}, in which we proposed a DRL-based framework for solving the MPP under demand uncertainty. The problem was formulated as a Markov decision process (MDP) that captures key combinatorial aspects, including demand uncertainty, vessel capacity, stability requirements, and the optimization of cargo revenue, hatch overstowage, and excess quay crane moves. To solve this formulation, we introduced an action policy with attention, also called attention model (AM) \cite{kool_attention_2019,vaswani_attention_2017}, with a feasibility projection layer, trained using DRL methods to produce adaptive and feasible solutions.

Building on this foundation, we propose AI2STOW, an end-to-end DRL policy for master stowage planning, which offers the following contributions:
\begin{itemize}
    \item \textbf{Extended MDP with Blocks}: We extend the original MDP with blocks to include paired block stowage patterns: an industrially relevant planning strategy often overlooked in existing stowage planning research \cite{van_twiller_literature_2024}. The extended implementation is released in the open-source repository\footnote{\href{https://github.com/OptimalPursuit/navigating_uncertainty_in_mpp}{GitHub Repository}}.
    \item \textbf{Action Mask to Enforce Paired Block Stowage}: We integrate an action-masking mechanism to enforce non-convex paired block stowage constraints in combination with projection layers to minimize convex feasibility violations in the DRL framework.
    \item \textbf{Efficient and Adaptive Solutions}: Experiments show that AI2STOW learns adaptive and feasible policies, outperforming baselines from stochastic programming and DRL in both objective quality and computational efficiency. The evaluation also includes a comparison of different projection layer configurations.
    \item \textbf{Decision Support for Realistic-Sized Instances}: AI2STOW can generalize well to larger problem instances, offering decision support for a realistic-sized vessel and operational planning horizons. These findings underscore the potential of DRL-based approaches in developing scalable algorithms for stowage planning.
\end{itemize}
\section{Background} \label{sec_drl_scale:background}
This section provides background knowledge on the domain of container vessel stowage planning and the preliminaries of reinforcement learning.

\subsection{Domain of Container Vessel Stowage Planning}\label{sec:domain}
Given the extensive nature of container stowage planning, we offer readers a general introduction to the domain. For a comprehensive exploration of the domain, we refer to \cite{jensen_container_2018}.

A liner shipping company operates a fleet of container vessels following fixed schedules on a closed-loop route. An apt analogy is that of a maritime bus line, transporting cargo between ports instead of people. Each complete journey to all ports is known as a voyage. The stowage planner typically receives information about the number of containers to be loaded at the next port several hours before the vessel's arrival. These voyages typically originate in loading regions (e.g., Asia) with a supply surplus and then sail towards discharge regions (e.g., Europe) characterized by a demand surplus. This dynamic causes a global container distribution imbalance, necessitating the redistribution of (empty) containers back to loading regions. The sailing leg from loading to discharge regions is named the \textit{fronthaul}, whereas the opposite is referred to as the \textit{backhaul}. 

The primary objective of stowage planning is to maximize cargo uptake based on the available vessel capacity, especially during the fronthaul. The secondary objective is to minimize operational costs and additional port fees, often associated with reducing the \textit{port stay} of the vessel. Modern vessels, with capacities exceeding 20,000 \textit{Twenty-Foot-Equivalent} (TEU) units, embark on voyages that typically include more than 10 ports. Due to its discrete and non-linear nature, the sheer size of this problem already poses a challenge.

Cargo typically comes in three sizes: 20 ft., 40 ft., and 45 ft., with dimensions of 8 ft. in width and 8 ft. 6 in. in height for \textit{standard containers} or 9 ft. 6 in. for \textit{highcubes}. The weight of containers ranges from 4 to 30 tonnes, depending on the load. Each container has an origin, known as the \textit{port of load} (POL), and a destination, referred to as the \textit{port of discharge} (POD).  Additionally, cargo types can vary, with \textit{Specials} including refrigerated containers \textit{reefers}, \textit{IMDG} containers carrying dangerous goods, and \textit{OOG} containers deviating from standard dimensions. For safety reasons, IMDGs are often segregated.

\begin{figure}[h!]
\centering
\includegraphics[scale=0.45]{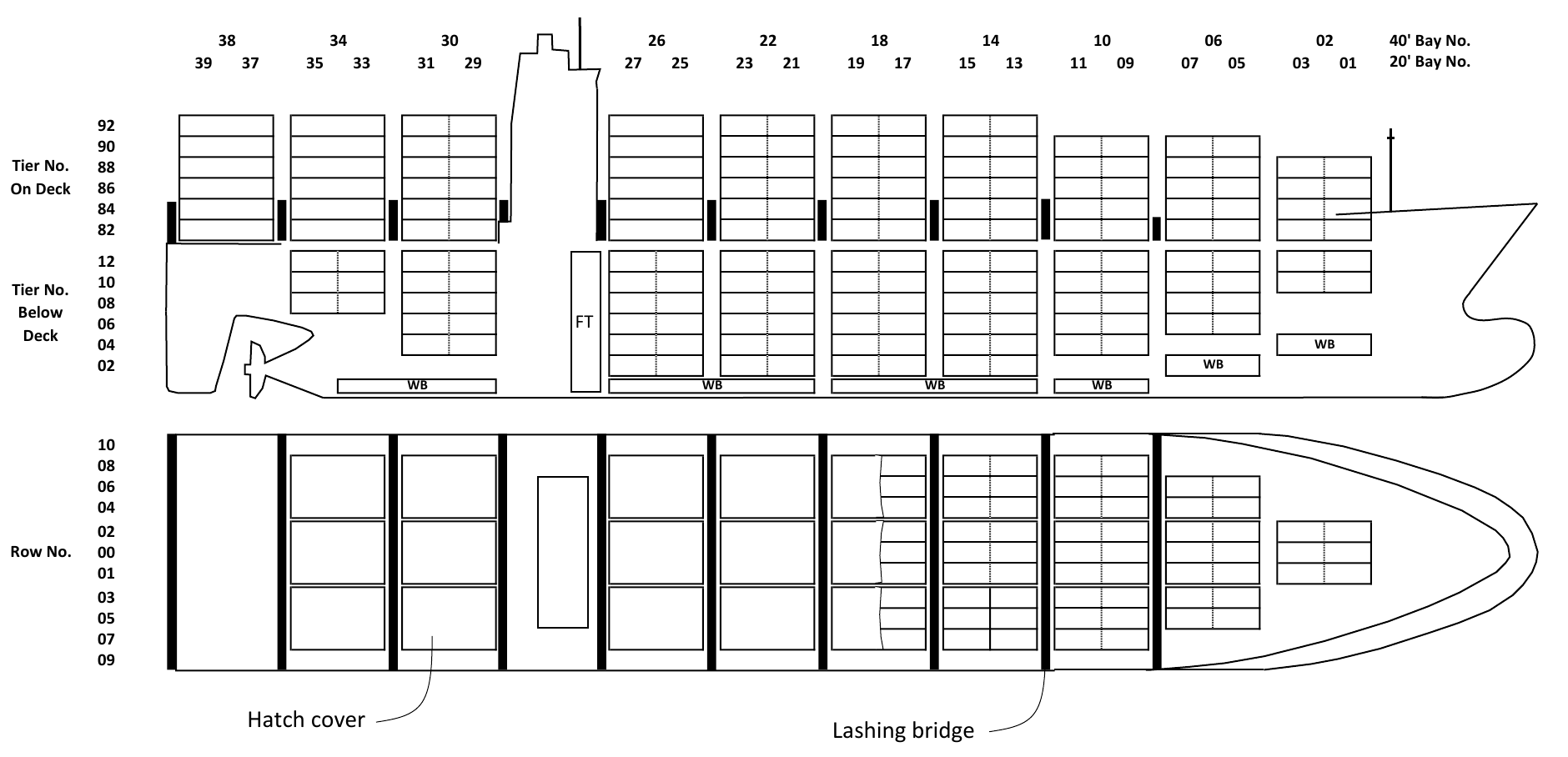}
\caption[Vessel side and top view]{Vessel side and top view \cite{van_twiller_literature_2024}}
\label{fig:vessel1}
\end{figure}

\begin{figure}[h!]
\centering
\includegraphics[scale=0.65]{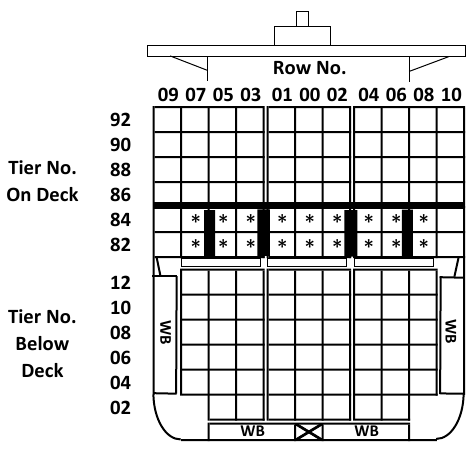}
\caption[Vessel front view]{Vessel front view \cite{van_twiller_literature_2024}}
\label{fig:vessel2}
\end{figure}

To ensure efficient transportation, container vessels adopt a cellular design, shown in Figures \ref{fig:vessel1} and \ref{fig:vessel2}. The vessel is compartmentalized into \textit{bays} (02 - 38) featuring \textit{rows} and \textit{tiers}, containing \textit{cells} capable of holding one 40-foot container or two 20-foot containers in \textit{slots}. Instead of front and back, we refer to the vessel's \textit{fore} and \textit{aft}, respectively. \textit{Stacks} are vertical arrangements of {cells}, whereas bays are horizontally separated by \textit{hatch covers} into \textit{locations} either \textit{on-deck} or \textit{below-deck}. Additionally, bays are vertically segmented into \textit{blocks} to combine all stacks above or below hatch covers.

To secure below-deck stacks (02 - 12), the vessel incorporates \textit{cell guides}, and the below-deck hold is sealed with \textit{hatch covers}. On-deck stacks (82 - 92) rest on hatch covers or the ship deck. \textit{Twist locks} bind stacked containers, while \textit{lashing rods} connect container corners to the deck or lashing bridges, enhancing stability. In addition to capacity limits, stacks have weight limits, with below-deck stacks also having height limits. Figure \ref{fig:vessel2} highlights cells (82-84) equipped with power plugs for reefers, denoted by the symbol $\ast$.

Vessels must comply with safety regulations to ensure seaworthiness, including the essential distribution of cargo with different weights along the vessel. It is imperative to avoid listing, prevent grounding, and have sufficient transverse stability, measured by \textit{trim}, \textit{draft}, and \textit{metacentric height} ({GM}). These factors, in turn, are influenced by the \textit{longitudinal} ({LCG}), \textit{vertical} ({VCG}), and \textit{transversal center of gravity} ({TCG}). Stress forces, including \textit{shear force}, \textit{bending moment}, and \textit{torsion moment}, also require careful consideration. Ultimately, each factor influences the vessel's capacity and must be kept at an acceptable level before embarking on a voyage.

In addition to factors affecting vessel capacity, considerations that impact port stay are also crucial. As previously mentioned, bays are separated by hatch covers into above- and below-deck locations. The cover must be removed to reach below-deck cargo in any bay, requiring \textit{quay cranes} to perform a \textit{restow} operation, moving on-deck containers that \textit{overstow} below-deck cargo. This phenomenon, known as \textit{hatch overstowage}, adds unnecessary moves to the port stay. Minimizing hatch overstowage is recognized as an NP-hard task, making it challenging to find an optimal solution in polynomial time \cite{tierney_complexity_2014}. Similarly, \textit{stack overstowage} can occur when stacks are not arranged in ascending order of PODs from top to bottom. 

Another factor influencing port stay is the \textit{makespan} of quay cranes operating in parallel. The makespan is defined by the \textit{long crane}, which is the crane with the longest operational time. Each port terminal commits to a certain number of moves within a specified time frame without guaranteeing a minimum number of cranes. This encourages stowage planners to distribute moves to ensure efficient quay crane operations. The \textit{crane intensity} (\textit{CI}) serves as an upper bound on the number of cranes that can work in parallel without idle time and is calculated as $\textit{CI} = \textit{TM}/\textit{WLC}$, where $\textit{TM}$ represents the total number of load and discharge moves, and $\textit{WLC}$ is the workload of the long crane in terms of moves. Stowage plans should maintain a CI greater than the average number of cranes assigned to the vessel.

Furthermore, planners use \textit{stowage patterns}, which are specific configurations designed to create robust plans ensuring high flexibility for cargo loaded in future ports and ensuring efficient crane operations. In \textit{block stowage}, for instance, two blocks separated by the hatch cover only contain containers with the same POD. \textit{Paired block stowage} (PBS) extends this notion to ensure the same POD in \textit{wing blocks} (rows 03-09 and 04-10 in Figure \ref{fig:vessel2}), while \textit{center blocks} (rows 00-02) can hold cargo with a different POD.

Until now, the focus has been on safety and cost reduction. However, revenue management and yield optimization are also important aspects of stowage planning. Generally, cargo can be categorized into two customer types, namely, \textit{long-term contracts} - high volume, lower variability and lower revenue - and \textit{spot market contracts} - lower volume, higher variability and higher revenue. Depending on the market dynamics, prices and demand will change over time. Due to the absence of no-show costs, both customer types have inherent uncertainty, as cargo may be booked onto vessels but will not arrive before the \textit{cargo cut-off time}. Typically, uncertainty reduces over time until this cut-off time, after which the demand is fully known. 

These considerations combine in a CO problem employed for constructing stowage plans \cite{jensen_container_2018}. Given the substantial scale of the challenge, this leads to a complex multi-port problem, for which  (near-)optimal solutions are still yet to be found \cite{van_twiller_literature_2024}.

\subsection{Preliminaries of Reinforcement Learning}
Reinforcement learning (RL) environments are commonly formalized as Markov decision processes (MDPs). An episodic discounted MDP is defined by \(\mathcal{M} =(S,X,\mathcal{T},\mathcal{R},T_\textit{epi}, \gamma)\), where ${s} \in S \subseteq \mathbb{R}^{n_s}$ represents observable states, and ${x} \in X \subseteq \mathbb{R}^{n_x}$ denotes actions executed by a policy. The state transition dynamics are stochastic by $\mathcal{T}: S \times X \rightarrow \Delta(S)$, where $\Delta(S)$ denotes the set of probability distributions over the state space $S$, or deterministic through $\mathcal{T}: S \times X \rightarrow S$. The reward function is defined by $\mathcal{R}: S\times X \rightarrow \mathbb{R}$, quantifying feedback for state-action pairs. The MDP has an episode length of $T_\textit{epi}$ and includes a discount factor $\gamma \in (0,1)$, which discounts future rewards to allow asymptotic convergence \cite{sutton_reinforcement_2018}.

Within the MDP, RL aims to derive an optimal policy that prescribes actions that maximize some objective. The policy $\pi({x}|{s})$ defines action probabilities given state observations. In deterministic scenarios, the policy simplifies to $\pi: S\rightarrow X$. Interaction between the agent and the environment unfolds sequentially in episodes: at each timestep $t$, the agent observes state ${s}_t$, executes action ${x}_t$, receives reward $r_{t+1}=\mathcal{R}({s}_t, {x}_t)$, and transitions to the next state are sampled by ${s}_{t+1}\sim\mathcal{T}({s}_{t+1}|{s}_t, {x}_t)$ with $\sim$ representing sampling. This cycle continues until it reaches a terminal state ${s}_{\text{term}}$ or timestep $T_{epi}$.

The objective of an optimal policy $\pi^*$ is to maximize expected discounted return or cumulative reward, hence, the objective given a policy $\pi$ is defined in Equation \eqref{eq:obj}.
\begin{align} \label{eq:obj}
\max_\pi J(\pi) = \mathbb{E}_\pi\left[ \sum_{t=0}^{\infty} \gamma^t r_{t} \right]
\end{align}
The agent's performance metric is the expected discounted return, as shown in Equation \eqref{eq:return}. Using the return $G_t$, we define performance for states via the value function $V({s}_t)$ in Equation~\eqref{eq:value}, and for state-action pairs through the Q-value function $Q({s}_t, {x}_t)$ in Equation~\eqref{eq:q_value}.
\begin{align}
G_t & = \sum_{h=0}^{T_{epi}} \gamma^h r_{t+h+1} \label{eq:return} \\
V^\pi({s}_t) & = \mathbb{E}_\pi \left[\sum_{h=0}^{T_{epi}} \gamma^k r_{t+h+1}\mid {s}_t\right], \label{eq:value} \\
Q^\pi({s}_t, {x}_t) & = \mathbb{E}_\pi \left[\sum_{h=0}^{T_{epi}} \gamma^k r_{t+h+1}\mid {s}_t, {x}_t\right]. \label{eq:q_value}
\end{align}

Solving MDPs is analogous to finding the an optimal policy $\pi^*$, which prescribes the best possible actions ${x}_t \sim \pi^*$ conditioned on the current state ${s}_t$. This optimal policy guides decision-making in each state ${s}_t$ to ensure the maximal cumulative rewards over time. Several groups of RL algorithms can be used to achieve this.

For instance, dynamic programming (DP) consists of model-based methods where $\mathcal{T}({s}'|{s},{x})$ and $\mathcal{R}({s},{x})$ are assumed to be known to estimate solution $V({s})$ by planning ahead as shown in Equation \eqref{for:value_bellman}. Although convergence to $\pi^*({x}|{s})$ is guaranteed, as per contraction mappings and the Banach fixed-point theorem \cite{puterman_markov_1994}. DP is intractable for large state and action spaces, and assumes perfect information that is often unrealistic \cite{sutton_reinforcement_2018}.
{
\begin{align}
    V({s}) & = \max_{x} \mathcal{R}({s},{x}) + \gamma\sum_{s'}\mathcal{T}({s}'|{s},{x})V({s}') \label{for:value_bellman}
\end{align}}

Instead, one could consider learning from experience using model-free methods. For example, Monte Carlo (MC) methods perform trial-and-error simulations of an episode and then update the value function $V({s}_t)$ based on the observed return $G_t$ from step $t$ until the end of the episode. Equation \eqref{for:monte_carlo} defines the update of $V({s}_t)$ at the end of episodes, where $\eta \in (0,1]$ is the learning rate of the task. This relatively straightforward concept can converge through sufficiently many simulations, as per the law of large numbers \cite{sutton_reinforcement_2018}. However, MC methods are susceptible to high variance, leading to unstable learning, particularly in environments with long episodes and high stochasticity \cite{sutton_reinforcement_2018}.
{ 
\begin{align} \label{for:monte_carlo}
V({s}_t) \leftarrow V({s}_t) + \eta \cdot (G_t - V({s}_t)) 
\end{align}}

Alternatively, one can consider temporal difference (TD) learning to obtain approximate functions with bootstrapping \cite{sutton_reinforcement_2018}. In this way, it is possible to update the value function $V({s}_t)$ at each timestep during an episode by looking ahead using knowledge of value estimates and future rewards. Equation \eqref{for:td_learning} defines the update of $V({s}_t)$, where $r_{t+1}$ is the immediate reward at timestep $t+1$. It should be noted that TD learning is proven to converge to $V^*({s}_t)$ under certain conditions, as per contraction mappings and the Banach fixed-point theorem \cite{sutton_reinforcement_2018}. Nonetheless, TD learning suffers from bias, potentially leading to systematic estimation errors. Additionally, its efficiency is sensitive to the initial values of $V({s}_t)$ \cite{sutton_reinforcement_2018}. 
{\begin{align} \label{for:td_learning} 
V({s}_t) \leftarrow V({s}_t) + \eta \cdot (r_{t+1} + \gamma \cdot V({s}_{t+1}) - V({s}_t))
\end{align}}

Another way is through DRL, which approximates solutions in large state and action spaces using parameterized functions (e.g, policy $\pi_\theta$ with parameter $\theta$). The unfortunate aspect is that convergence guarantees of Equations~\eqref{for:value_bellman}, \eqref{for:monte_carlo}, and \eqref{for:td_learning}
 will be lost, as neural networks as function approximators break the required assumptions of convergence guarantees. Instead, the policy gradient theorem states that the derivative of the objective $\nabla J(\theta)$ with parameter $\theta$ is proportional to the terms in Equation \eqref{for:PG_theorem}, where $\psi_\pi({s}_t)$ is the on-policy probability of state ${s}_t$. This allows learning \(\pi_\theta\) by updating \(\theta\) with respect to the objective \(J(\theta)\) \cite{sutton_reinforcement_2018}. Generally, the objective is to maximize the expected return under a parameterized policy as $J(\theta)=\mathbb{E}_{\pi_\theta} [\sum_{t=0}^{\infty} \gamma^t r_{t}] $. Due to large policy updates, gradient-based methods can suffer from variance, hindering stability and efficiency of learning \cite{sutton_reinforcement_2018}. To mitigate this variance, a baseline (e.g., $V(s_t)$) is subtracted from the return, preserving the unbiasedness of the gradient while improving the stability and efficiency of learning.
{ \begin{align} \label{for:PG_theorem}
    \nabla J(\theta) & \propto \mathbb{E}_{\pi}\left[\sum_{{s}_t} \psi^{\pi}({s}_t) \sum_{{x}_t} \ln \nabla \pi_\theta({x}_t|{s}_t) Q^\pi({s}_t,{x}_t) \right] 
\end{align}}

Furthermore, actor-critic methods can mitigate the drawbacks of learning from experience and policy gradients by combining both techniques \cite{sutton_reinforcement_2018}. The critic estimates the expected return \( V(s_t) \), serving as a baseline to reduce variance, while the actor adjusts the policy \( \pi(x_t \mid s_t) \) based on this evaluation of the \((s_t, x_t)\)-pair. The characteristics of actor-critic methods depend on the specific algorithms used but generally are relatively sample-efficient with reliable performance \cite{schulman_proximal_2017}.
Hence, we favour actor-critic methods due to their ability to leverage the strengths and mitigate the drawbacks of policy gradients and TD-learning. 
\section{Related Work} \label{sec_drl_scale:related_work}
This section covers relevant contributions in container stowage optimization, as well as the use of machine learning in the context of CO.

\subsection{Container Stowage Planning} 
The complexity of container vessel stowage planning arises from its inherent multi-port nature. Consequently, algorithms are required to balance a myriad of combinatorial aspects at each port. The field can be broadly categorized into single-port work aimed at creating light-weight operational stowage plans (e.g., \cite{ambrosino_stowing_2004, ambrosino_experimental_2010, delgado_placement_2012, larsen_heuristic_2021}), and multi-port initiatives dedicated to creating realistic stowage plans (e.g., \cite{botter_stowage_1992, avriel_stowage_1998, wilson_container_2000, pacino_fast_2011, roberti_decomposition_2018, parreno-torres_solving_2021}). As analyzed in \cite{van_twiller_literature_2024}, a definitive and optimal solution to either problem is yet to be found.

\begin{figure}[h!]
\centering
\includegraphics[scale=0.4]{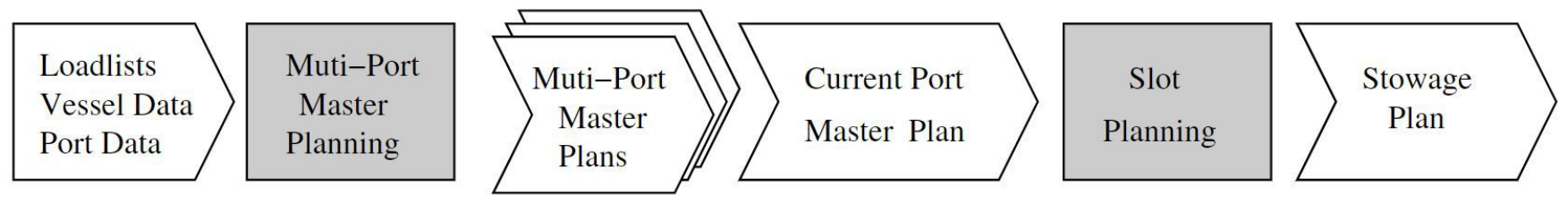}
\caption[Hierarchical decomposition of stowage planning]{Hierarchical decomposition of stowage planning \cite{pacino_fast_2011}}
\label{fig_drl_scale:decomp}
\end{figure}

To address the inherent complexity of stowage planning, a recommended strategy involves hierarchical decomposition, which systematically divides the problem into sequential subproblems. One widely recognized version, shown in \Cref{fig_drl_scale:decomp}, decomposes the problem into the master planning problem (MPP) and the slot planning problem (SPP) (e.g., \cite{wilson_container_2000, pacino_fast_2011}). During the MPP, groups of containers are allocated to locations on the vessel while aiming to satisfy global objectives and constraints (e.g., \cite{pacino_lns_2013, bilican_mathematical_2020, chao_minimizing_2021}). Subsequently, the focus shifts to the SPP, where individual containers are assigned to specific slots within the designated locations that aim to satisfy local objectives and constraints (e.g., \cite{pacino_3-phase_2010, kebedow_including_2019, korach_matheuristics_2020}). For a comprehensive understanding of hierarchical decomposition, readers are referred to \cite{jensen_container_2018}. Furthermore, most stowage planning research assumes deterministic cargo demand, whereas the more realistic stochastic variant has received limited attention, with only a single study addressing it \cite{christensen_rolling_2019}. Recently, heuristic frameworks have gained traction as alternatives to hierarchical decomposition, offering more ways to handle the complexity of stowage planning (e.g., \cite{pacino_crane_2018, larsen_heuristic_2021}). 

Regardless, several solution methods have been proposed to solve different stowage planning problems, for instance, exact methods  (e.g., \cite{roberti_decomposition_2018,zhu_integer_2020}), greedy heuristics (e.g., \cite{avriel_stowage_1998,ding_stowage_2015}), population-based  (e.g., \cite{dubrovsky_genetic_2002,chang_solving_2022}) or neighborhood-based metaheuristics (e.g., \cite{ambrosino_experimental_2010,pacino_crane_2018}), matheuristics (e.g., \cite{korach_matheuristics_2020,parreno-torres_solving_2021}), tree-based methods (e.g., \cite{azevedo_solving_2014}), or hybrid frameworks  (e.g., \cite{wilson_container_2000,bilican_mathematical_2020}).

To the best of our knowledge, the number of contributions related to RL in stowage planning is limited. Most contributions have focused on the placement of individual containers in single-port problems that omit key combinatorial aspects by methods such as deep Q-learning \cite{shen_deep_2017} and Monte Carlo tree search \cite{zhao_container_2018}. Given the scale of modern vessels, applying RL at the container level would require placing up to 20,000 containers per voyage, obtaining very long episodes that complicate learning. Therefore, focusing on the MPP, which abstracts the problem to a higher level, may offer a more practical and scalable alternative. To date, the only research applying RL to the MPP has been conducted by us \cite{van_twiller_towards_2023, van_twiller_navigating_2025}. In these studies, we model the MPP as an MDP and use a DRL approach to derive solutions. While promising, this work requires further improvements in representativeness. To address this, we aim to extend the MPP formulation by incorporating paired block stowage patterns (PBS) and accounting for demand uncertainty while ensuring scalability to industrial vessel sizes and realistic voyage lengths. Despite the growing algorithmic diversity, the search continues for scalable, generalizable methods capable of solving representative stowage planning problems.

\subsection{Stochastic Programming}
A classical approach for optimization under uncertainty is stochastic programming, where uncertainty is modeled explicitly through a discrete set of scenarios that evolve over stages of time \cite{birge_introduction_2011}. This results in a scenario tree representation, whose size grows exponentially with the number of stages, $\mathcal{O}(b^T)$, where $b$ is the branching factor and $T$ is the number of stages. As a result, multi-stage formulations become computationally challenging beyond a few stages.

Numerous techniques have been developed to address these challenges. Scenario reduction methods \cite{dupacova_scenario_2003, watanabe_scenario_2009} aim to approximate the full scenario tree by selecting a representative subset of scenarios that preserves key probabilistic characteristics. Decomposition techniques, including Benders decomposition and variants \cite{rahmaniani_benders_2017} or Lagrangian relaxation \cite{caroe_dual_1999}, divide the original problem into smaller, more manageable subproblems that are solved iteratively while coordinating shared constraints. The progressive hedging algorithm, introduced by \cite{rockafellar_scenarios_1991} and extended in later work \cite{boland_combining_2018}, solves scenario subproblems independently while enforcing non-anticipativity through augmented Lagrangian penalties. Stochastic dual dynamic programming\cite{pereira_multi-stage_1991,shapiro_analysis_2011} approximates value functions in multi-stage problems by iteratively generating Benders cuts through backward passes and simulating decision paths in forward passes. The sample average approximation method \cite{shapiro_monte_2003,chen_sample_2022} replaces expectations with sample averages over finite scenario sets, enabling tractable optimization and offering convergence guarantees as the sample size increases.

Applications span domains such as supply chain optimization \cite{ahmed_convexity_2006,thevenin_stochastic_2022}, unit commitment in power systems \cite{papavasiliou_stochastic_2012}, and telecommunications network design \cite{atamturk_two-stage_2007}. Recent trends include data-driven stochastic programming, which leverages historical data or machine learning for scenario generation and decision policy learning \cite{ban_big_2019, bertsimas_data-driven_2018}. Another emerging direction is distributionally robust optimization, where solutions are sought under worst-case distributions within an ambiguity set \cite{mohajerin_esfahani_data-driven_2018, gao_distributionally_2023}.

Despite substantial progress, scaling stochastic programming to high-dimensional, multi-stage, and real-time decision-making environments remains challenging, motivating further research into hybrid models, learning-based approximations, and scalable decomposition frameworks.

\subsection{Machine Learning for Optimization Problems}
In recent years, ML has proven to be an effective tool for solving CO problems, at times outperforming or enriching conventional solution methods \cite{bengio_machine_2021}. Specifically, DRL has emerged as a promising method to deal with challenging CO problems \cite{mazyavkina_reinforcement_2021}. The usage of ML can be classified into learning categories, each of which will be discussed in the following paragraphs.

End-to-end learning leverages ML to directly output solutions for input problem instances,  thereby circumventing the need for handcrafted heuristics or manually designed search procedures. One particular example can be found in the chip design process, which can also be decomposed hierarchically, with chip floor planning being one sequential subproblem. In this context, a graph convolutional encoder with actor-critic decoders generates an approximate chip floor plan that dictates chip quality, leading to a significant acceleration of the overall process with chips of equal or superior quality compared to human-made designs \cite{mirhoseini_graph_2021}. Additionally, actor-critic methods have demonstrated excellent performance on various control tasks \cite{schulman_proximal_2017,mirhoseini_graph_2021,openai_gpt-4_2023}. Other approaches to solving well-known combinatorial optimization problems include pointer networks that integrate recurrent neural networks (RNNs) with attention mechanisms \cite{vinyals_pointer_2015}, a graph-based encoder combined with an RNN decoder using multi-head attention and optimized via an actor-critic REINFORCE algorithm \cite{nazari_reinforcement_2018}, and graph attention models trained with REINFORCE using a greedy rollout baseline \cite{kool_attention_2019} and multiple rollouts \cite{kwon_pomo_2020}.

Within learning solution heuristics, a key challenge is ensuring feasibility in complex and explicit action spaces, particularly when feasible regions are dynamic and state-dependent. Various deep learning methods have addressed constraint learning through techniques such as backpropagation over (in)equality completions \cite{donti_dc3_2021}, differentiable projection layers that map interior points to boundary regions \cite{li_learning_2023}, convex programming layers \cite{agrawal_differentiable_2019}, and problem-specific repair mechanisms \cite{chen_end--end_2024}. In addition, safe reinforcement learning has approached feasibility through constrained MDPs with primal-dual techniques \cite{ding_natural_2020}, soft barrier functions \cite{wang_enforcing_2023}, and safety shields \cite{alshiekh_safe_2018}.

Even though end-to-end learning can be advantageous, learning to configure algorithmic components can also be effective in guiding the search process. For example, efficient active search updates the weights of solution-constructing models via DRL for problems such as capacitated vehicle routing and job shop scheduling \cite{hottung_efficient_2022}. Another approach uses an attention model trained with REINFORCE to reconstruct neighborhoods within a large neighborhood search (LNS) framework, outperforming methods with handcrafted heuristics on capacitated and split delivery vehicle routing problems \cite{hottung_neural_2019}.

Another way is to use ML and CO algorithms in parallel. A literature survey explores traditional and learning methods for variable and node selection in branch-and-bound frameworks \cite{lodi_learning_2017}. For example, deep learning can learn solution strategies and lower bounds by analyzing existing (near-)optimal solutions to instances, which are integrated into a tree search framework to assist branching and pruning decisions \cite{hottung_deep_2020}.  
\section{Problem Formulation} \label{sec_drl_scale:problem}
Our problem formulation models the key combinatorial aspects of the MPP, aiming to generate an approximate plan that meets global objectives and constraints under demand uncertainty. These include maximizing cargo revenue, minimizing hatch-overstowage and excess crane move costs, applying valid PBS patterns, and ensuring acceptable values for long crane length, metacentric height, and trim. Like the early stages of stowage planning, where primary objectives take precedence, this problem focuses on global objectives and constraints during the voyage. It intentionally abstracts away individual containers and slots, allowing the interchangeable placement of similar containers in various slots. The solution to the MPP guides the subsequent container placement, considering specific local objectives and constraints in the SPP.

\Cref{tab_drl_scale:sets} defines the sets used in this MPP. Let us define a voyage denoted by the ordered set of ports $P$ with $N_P$ being the last port. Given two ports in $P$, we can derive a subset of the voyage ${P}_\textit{start}^\textit{end}$. Let $\textit{TR}$ be the set of transport pairs of all possible POLs and PODs. Given port $p$, various subsets of $\textit{TR}$ can be defined, namely, transports on board before continuing the voyage ($\textit{TR}^\textit{OB}(p)$), the transports remaining on board (ROB) after discharging but before loading ($\textit{TR}^\textit{ROB}(p)$), and all transports that are either loaded or discharged ($\textit{TR}^\textit{M}(p)$).  Furthermore, containers are classified into distinct cargo classes in the set $K$, which include container size (1 or 2 TEU), weight classes (light, medium and heavy) and customer types (spot market or long-term contracts). Consider that vessel locations are defined by bays, decks and paired blocks. Let $B$ be the ordered set of bays with $N_B$ being the number of bays, indicating the longitudinal position extending from fore to aft. For any adjacent pair of bays, we define the set of bay pairs $B'$. Let $D$ represent the set of decks, differentiating between on-deck $d_{o}$ and in the hold $d_{h}$, thereby determining the vertical positioning. Let $\textit{BL}$ be the set of paired blocks with $N_\textit{BL}$ being the number of paired blocks in a bay. If \( N_\textit{BL} = 1 \), there is only a single block. For \( N_\textit{BL} > 1 \), the first paired blocks in a bay correspond to the wing blocks, while the remaining blocks are center blocks. Note that modern vessels have \( N_\textit{BL} \leq 3 \).

\begin{table}[h!]
\centering
\small
\caption{Sets of the MPP}
\label{tab_drl_scale:sets}
\begin{tabular}{ll}
\toprule
Ports &  $p \in P=\{1,2,\ldots, N_P\}$ \\
Port range & $p \in {P}_\textit{start}^\textit{end} = \{p \in {P}\!\mid\! \textit{start}\leq p \leq \textit{end}\}$ \\  
Transport pairs &  $\textit{tr} \in \textit{TR}=\{(i,j) \in P^2\mid i<j\}$ \\
Onboard transp. &  $\textit{tr} \in \textit{TR}^{\mathit{OB}}(p)=\{(i,j) \in P^2\mid i\leq p, j > p\}$ \\
ROB transports &  $\textit{tr} \in \textit{TR}^{\mathit{ROB}}(p)=\{(i,j) \in P^2\mid i<p, j>p\}$ \\
Discharge moves & $\textit{tr} \in \textit{TR}^{-}(p) = \{(i,p) \in P^2\mid i < p\}\;$ \\
Load moves & $\textit{tr} \in \textit{TR}^{+}(p) = \; \{(p,j) \in P^2\mid j > p\}$ \\
Port moves & $\textit{tr} \in \textit{TR}^{M}(p) = \textit{TR}^{+}(p) \cup \; \textit{TR}^{-}(p)$ \\
Cargo classes &  $k \in {K} = {\scriptsize \{\textit{20ft}, \textit{40ft}\} \times \{\textit{Light}, \textit{Medium}, \textit{Heavy}\} \times \{\textit{Spot},\textit{Long}\}}$ \\
Bays &  $b \in B =  \{1,2,\ldots, N_B\}$ \\
Adjacent bays & $b' \in B' = \{(1,2),(2,3),\ldots,(N_B-1, N_B)\}$ \\
Decks &  $d \in D = \{d_{o}, d_{h}\}$ \\
Paired blocks &  $\textit{bl} \in \textit{BL} = \{1, 2, \dots, N_\textit{BL}\}$ \\
\bottomrule
\end{tabular}
\end{table}

\Cref{tab_drl_scale:params} defines the MPP parameters. The transport matrix \( q \) represents realized demand in containers, \(\textit{rev}\) is the revenue per transport and container type, \(\textit{ct}^\textit{ho} \) refers to the hatch overstowage cost parameter, and \(\textit{ct}^\textit{cm} \) is the excess crane moves cost parameter. Additionally, \( c \) denotes vessel location capacity in TEU, while container sizes and cargo weights are given by \( \textit{teu} \) and \( w \), respectively. The longitudinal distance \( \mathit{ld} \) is measured from fore to aft, and the vertical distance \( \mathit{vd} \) is measured from keel upwards, with the center normalized to 1. Stability constraints include the acceptable bounds for the longitudinal and vertical center of gravity (\( \underline{\textit{lcg}}, \overline{\textit{lcg}} \) and \( \underline{\textit{vcg}}, \overline{\textit{vcg}} \)), as well as the permissible error in excess crane moves (\( \delta^\textit{cm} \)). We also define big $M$ as a large constant to enforce logical conditions.

\begin{table}[h!]
\centering
\small
\caption{Parameters of the MPP}
\label{tab_drl_scale:params}
\begin{tabular}{ll}
\toprule
Transport matrix of realized demand &  $q \in \mathbb{Z}^{|\textit{TR}|\times|K|}_{\geq 0}$\\ 
Revenue per container &  $\textit{rev} \in \mathbb{R}^{|\textit{TR}|\times|K|}_{> 0} $ \\
Cost per overstowed container &  $\textit{ct}^\textit{ho} \in \mathbb{R}_{> 0} $ \\
Cost per excess crane move &  $\textit{ct}^\textit{cm} \in \mathbb{R}_{> 0}$ \\
Location capacity in TEU &  $c \in \mathbb{Z}^{|B|\times|D| \times |\textit{BL}|}_{\geq 0} $ \\
Container size in TEU &  $\textit{teu} \in \mathbb{Z}^{|K|}_{\geq 0} $ \\
Cargo weight class in tonnes & $w \in \mathbb{R}^{|K|}_{\geq 0}$ \\ 
Longitudinal distance of bays & $\mathit{ld}_b = \frac{2b - 1}{|B|} \; \forall b \in B$\\
Vertical distance of deck & $\mathit{vd}_d = \frac{2d - 1}{|D|} \; \forall d \in D$\\
Acceptable LCG bounds &  $\underline{\textit{lcg}},\overline{\textit{lcg}} \in \mathbb{R}_{\geq 0}$ \\
Acceptable VCG bounds &  $\underline{\textit{vcg}},\overline{\textit{vcg}} \in \mathbb{R}_{\geq 0}$ \\
Acceptable error excess crane moves & $\delta^\textit{cm} \in \mathbb{R}_{\geq 0}$ \\ 
Big M & $M$ \\ 
\bottomrule
\end{tabular}
\end{table}

Let us specify some parameters based on the following definitions. The realized demand values \( q_{(i,j),k} \) are sampled from a generator \( \mathcal{G} \) parameterized by \( \vartheta_{i,j,k} \):
\begin{align}
q_{(i,j),k} \sim \mathcal{G}(\vartheta_{i,j,k}) \quad \forall (i,j) \in \textit{TR},\; k \in K,
\end{align}

The revenue of containers with type $k$ on transport $(i,j)$ is defined by the length of the transport $(j-i)$ plus some standard revenue $\textit{SR}$. Note that long-term contracts have their transport revenue $(j-i)$ reduced by the parameter $\textit{LR}$. 
\begin{align}\textit{rev}_{(i, j), k} =
\begin{cases} 
(j - i)(1 - \textit{LR}) \!+\! \textit{SR}, & \text{if } \textit{Long} \in k \\
(j - i) \! + \! \textit{SR}, & \text{if } \textit{Spot} \in k
\end{cases} \; \forall (i,j) \in \textit{TR}, k \in K.
\end{align}

The TEU per container depends on type $k$ as:
\begin{align}
\textit{teu}_k =
\begin{cases} 
1, & \text{if } \textit{20ft} \in k \\
2, & \text{if } \textit{40ft} \in k
\end{cases} \; \forall k \in K.
\end{align}

Similarly, container weight is determined by the weight class:
\begin{align}
\textit{w}_k =
\begin{cases} 
1, & \text{if } \textit{Light} \in k \\
2, & \text{if } \textit{Medium} \in k \\
3, & \text{if } \textit{Heavy} \in k
\end{cases} \; \forall k \in K.
\end{align}

To simplify the problem, we make the following assumptions: 
\begin{itemize}
\setlength{\itemsep}{0pt}
\item Each voyage starts with an empty vessel, which is box-shaped with equal capacity in each bay.
\item The voyages include ports with balanced loading and discharging operations rather than having ports that specialize in either operation.
\item Demand at the current port is fully observable and deterministic, whereas demand at subsequent ports is uncertain. 
\item The uncertainty is exogenous, arising from external factors beyond the system’s control and independent of model decisions.
\item Special container types (e.g., reefers and IMDG) are not considered. We assume that each location has sufficient capacity for reefers, while other specials with local constraints can be considered in an SPP phase.
\item Loading and discharge times are equal for all ports and cargo classes. 
\item Revenue of cargo uptake \(\textit{rev}_{\textit{tr},k}\) is larger than costs of overstowage  \(\textit{ct}^\textit{ho}\) or excess crane moves \(\textit{ct}^\textit{cm}\).
\end{itemize}

Stochastic programming approaches typically represent uncertainty through an explicit scenario tree. This is a directed tree \(T_\textit{ST} = (V_\textit{ST}, E_\textit{ST})\), where \( V_\textit{ST} \) is the set of nodes, each corresponding to a decision or uncertainty realization at a given stage. \( E_\textit{ST} \subseteq V_\textit{ST} \times V_\textit{ST} \) is the set of directed edges representing transitions between nodes over time. The tree $T_\textit{ST}$ consists of:  
\begin{enumerate}
    \item A root node \( v_1 \in V_\textit{ST} \), representing the initial state at the first port. 
    \item Stages \( p = 1, 2, \dots, N_{P}-1 \), where each node \( v \) belongs to a stage \( p(v) \). We denote stages by $p$, as stages are equivalent to ports in a voyage.   
    \item Branching structure, where each node has $S_\textit{ST}$ child nodes representing possible future realizations.
    \item A measure \( P_\textit{ST}: V_\textit{ST} \to [0,1] \) assigning probabilities to nodes, ensuring:     \begin{align}
   \sum_{v' \in \textit{child}(v)} \mathbb{P}(v') =  \mathbb{P}(v), \quad \forall v \in V_\textit{ST}.
   \end{align}  
   \item Scenario paths $\phi \in \mathcal{Z}$, which are root-to-leaf paths representing possible realizations of uncertainty over time. The number of scenario paths grows exponentially as $|\mathcal{Z}| = S_\textit{ST}^{N_P - 1}$ . 
\end{enumerate}

Based on the scenario tree, we define the MPP under demand uncertainty as a multi-stage stochastic MIP. \Cref{tab_drl_scale:vars} defines all decision variables under scenario path $\phi \in \mathcal{Z}$ provided by the scenario tree. Vessel utilization $\tilde{u}$ assigns cargo of type $k$ and transport $\textit{tr} = (i,j)$, with \textit{POL} $i$ and \textit{POD} $j$, to locations defined by bay $b$, deck $d$ and block $\textit{bl}$. Hatch overstowage $\tilde{\textit{ho}}$ represents the amount of containers that need to be restowed in bay $b$ and block $\textit{bl}$ at port $p$, while excess crane moves $\tilde{\textit{cm}}$ defines the amount of additional crane moves outside moves allocated to the vessel by a terminal at port $p$. Hatch movement $\tilde{\textit{hm}}$ indicates the need to access below-deck locations at bay $b$, block $\textit{bl}$ and port $p$, whereas the destination indicator $\tilde{\textit{di}}$ signifies that POD $j$ is assigned to bay $b$ and block $\textit{bl}$ at port $p$. 

\begin{table}[h!]
\centering
\small
\caption{Decision variables of the MPP}
\label{tab_drl_scale:vars}
\begin{tabular}{ll}
\toprule
Vessel utilization  &  $\Tilde{u}^{b,d,\textit{bl},\phi}_{\textit{tr},k} \in \mathbb{Z}_{\geq0} \; \forall b \in B, d \in D, \textit{bl} \in \textit{BL}, \textit{tr} \in \textit{TR}, k \in K, \phi \in \mathcal{Z}$\\
Hatch overstowage & $\Tilde{\textit{ho}}^{b,\textit{bl},\phi}_{p} \in \mathbb{Z}_{\geq0} \; \forall b \in B, \textit{bl} \in \textit{BL}, p \in P, \phi \in \mathcal{Z} $ \\
Excess crane moves & $\Tilde{\textit{cm}}^{\phi}_{p} \in \mathbb{Z}_{\geq0} \;  \forall p \in P, \phi \in \mathcal{Z}$ \\
Hatch movement & $\Tilde{\textit{hm}}^{b,\textit{bl},\phi}_{p} \in \{0,1\} \; \forall b \in B, \textit{bl} \in \textit{BL}, p \in P, \phi \in \mathcal{Z} $ \\
Destination indicator & $\Tilde{\textit{{di}}}^{b,\textit{bl},\phi}_{p,j} \in \{0,1\} \; \forall b \in B, \textit{bl} \in \textit{BL}, p \in P, j \in P^{N_P}_p, \phi \in \mathcal{Z} $ \\
\bottomrule
\end{tabular}
\end{table}

The objective function \eqref{for_drl_scale:MIP_obj} maximizes the expected revenue with parameter \(\textit{rev}_{\textit{tr},k} \in \mathbb{R}_{>0}\) and minimizes the expected hatch-overstowage with parameter \( \textit{ct}^\textit{ho} \in \mathbb{R}_{>0} \) and the expected crane moves costs with parameter \( \textit{ct}^\textit{cm} \in \mathbb{R}_{>0} \) over scenario paths $\phi \in \mathcal{Z}$ with probability $\mathbb{P}_\phi$. We assume equal probability for each scenario path.

Constraint \eqref{for_drl_scale:MIP_demand} limits onboard utilization to the cargo demand $q$, whereas Constraint \eqref{for_drl_scale:MIP_capacity} limits each vessel location to the TEU capacity $c$ for each location. Constraint \eqref{for_drl_scale:pods} links utilization to the destination indicator, while Constraint \eqref{for_drl_scale:max_pods} only allows a single POD for each bay and block, thereby enforcing PBS. In Constraint \eqref{for_drl_scale:hatch}, we indicate that hatches need to be opened if below deck cargo needs to be loaded or discharged. Based on these movements, Constraint \eqref{for_drl_scale:hatch_restow} models the amount of hatch overstowage in containers. Subsequently, we compute the target of crane moves $\overline{z}$ in Constraint \eqref{for_drl_scale:z_upper}, after which Constraint \eqref{for_drl_scale:long_crane} computes the excess number of crane moves $\Tilde{\textit{cm}}$. Additionally, we model the longitudinal and vertical stability in Constraints \eqref{for_drl_scale:lm} until \eqref{for_drl_scale:VS1}. First, we compute the longitudinal moment, vertical moment and total weight in Constraints \eqref{for_drl_scale:lm}, \eqref{for_drl_scale:vm} and \eqref{for_drl_scale:TW}, respectively. Second, Constraint \eqref{for_drl_scale:LS1} bounds \textit{lcg} between $\underline{\textit{lcg}}$ and $\overline{\textit{lcg}}$. Third,  Constraint \eqref{for_drl_scale:VS1} bounds \textit{vcg} between $\underline{\textit{vcg}}$ and $\overline{\textit{vcg}}$. Finally, we include non-anticipation in Constraint \eqref{for_drl_scale:nonanti} to prevent leveraging future demand realizations.
\begin{align}
\text{max } & 
\sum_{\phi \in \mathcal{Z}} \mathbb{P}_\phi \sum_{p\in {P}} \sum_{b \in{B}}\sum_{d\in {D}} \sum_{\textit{bl} \in{\textit{BL}}}\sum_{k\in {K}}\sum_{\textit{tr} \in \textit{TR}^\textit{+}(p)} \textit{rev}_{\textit{tr},k} \Tilde{u}^{b,d,\textit{bl},\phi}_{\textit{tr},k} \qquad\qquad \nonumber \\
& \qquad\qquad\qquad\qquad\qquad\qquad\qquad - \textit{ct}^\textit{ho} \Tilde{\textit{ho}}^{b,\textit{bl},\phi}_{p} - \textit{ct}^\textit{cm} \Tilde{\textit{cm}}^{\phi}_{p}\label{for_drl_scale:MIP_obj} \\
    \text{s.t. } & 
    \sum_{b \in {B}} \sum_{d \in {D}} \sum_{\textit{bl} \in \textit{BL}} \Tilde{u}^{b,d,\textit{bl},\phi}_{\textit{tr},k} \leq q^{\phi}_{\textit{tr},k} \nonumber\\
    & \qquad\qquad\qquad \forall p \in {P}, \textit{tr} \in \textit{TR}^\textit{OB}(p), k \in {K}, \phi \in \mathcal{Z} \label{for_drl_scale:MIP_demand} \\
    & \sum_{k \in {K}} \sum_{\textit{tr} \in \textit{TR}^\textit{OB}(p)} \textit{teu}_{k} \Tilde{u}^{b,d,\textit{bl},\phi}_{\textit{tr},k} \leq c_{b,d, \textit{bl}} \nonumber\\
    & \qquad\qquad\qquad \forall p \in {P}, b \in{B}, d \in {D}, \textit{bl} \in {\textit{BL}}, \phi \in \mathcal{Z} \label{for_drl_scale:MIP_capacity} \\
    & \sum_{k \in {K}} \Tilde{u}^{b,d,\textit{bl},\phi}_{(p,j),k} \leq {M} \Tilde{\textit{{di}}}^{b,\textit{bl},\phi}_{p,j} \nonumber\\
    & \qquad\qquad\qquad \forall p \in {P}, j \in {P}_p^{N_P},  b \in{B}, d \in{D},  \textit{bl} \in \textit{BL}, \phi \in \mathcal{Z} \label{for_drl_scale:pods} \\
    & \sum_{j \in {P}_p^{N_P}} \Tilde{\textit{{di}}}^{b,\textit{bl},\phi}_{p,j} \leq 1 \nonumber\\
    & \qquad\qquad\qquad \forall p \in {P},  b \in{B}, \textit{bl} \in \textit{BL}, \phi \in \mathcal{Z} \label{for_drl_scale:max_pods} \\
    & \sum_{k \in {K}} \sum_{\textit{tr} \in \textit{TR}^{M}(p)} \Tilde{u}^{b,d_h,\textit{bl},\phi}_{\textit{tr},k} \leq {M} \Tilde{\textit{hm}}^{b,\textit{bl},\phi}_{p} \nonumber\\
    & \qquad\qquad\qquad \forall p \in {P}, b \in{B}, \textit{bl} \in \textit{BL},\phi \in \mathcal{Z} \label{for_drl_scale:hatch} \\
    & \sum_{k \in {K}} \sum_{\textit{tr} \in \textit{TR}^\textit{ROB}(p)} \Tilde{u}^{b,d_o,\textit{bl},\phi}_{\textit{tr},k} - M(1 - \Tilde{\textit{hm}}^{b,\textit{bl},\phi}_{p}) \leq \Tilde{\textit{ho}}^{b,\textit{bl},\phi}_{p} \nonumber\\
    & \qquad\qquad\qquad \forall p \in {P}, b \in{B}, \textit{bl} \in \textit{BL}, \phi \in \mathcal{Z} \label{for_drl_scale:hatch_restow} \\
    & {\overline{z}}^{\phi}_{p} = (1 + \delta^\textit{cm}) \frac{2}{|{B}|} \sum_{\textit{tr} \in \textit{TR}^M(p)} \sum_{k \in {K}} q^{\phi}_{\textit{tr},k} \nonumber\\
    & \qquad\qquad\qquad \forall p \in {P}, \phi \in \mathcal{Z} \label{for_drl_scale:z_upper} \\
    & \sum_{b \in{b'}} \sum_{d \in {D}} \sum_{\textit{bl} \in \textit{BL}} \sum_{k \in {K}} \sum_{\textit{tr} \in \textit{TR}^{M}(p)} \Tilde{u}^{b,d,\textit{bl},\phi}_{\textit{tr},k} - {\overline{z}}^{\phi}_{p} \leq \Tilde{\textit{cm}}^{\phi}_{p} \nonumber\\
    & \qquad\qquad\qquad \forall p \in {P}, {b'} \in {B}', \phi \in \mathcal{Z} \label{for_drl_scale:long_crane} \\
    & \textit{tw}^{\phi}_{p} = \sum_{k \in {K}} w_k \sum_{\textit{tr} \in \textit{TR}^\textit{OB}(p)} \sum_{\textit{bl} \in \textit{BL}} \sum_{d \in {D}} \sum_{b \in{B}} \Tilde{u}^{b,d,\textit{bl},\phi}_{\textit{tr},k} \nonumber\\
    & \qquad\qquad\qquad \forall p \in {P}, \phi \in \mathcal{Z} \label{for_drl_scale:TW} \\
    & \textit{lm}^{\phi}_{p} = \sum_{b \in{B}} \textit{ld}_b \sum_{k \in {K}} w_k \sum_{\textit{tr} \in \textit{TR}^\textit{OB}(p)} \sum_{\textit{bl} \in \textit{BL}} \sum_{d \in {D}} \Tilde{u}^{b,d,\textit{bl},\phi}_{\textit{tr},k} \nonumber\\
    & \qquad\qquad\qquad \forall p \in {P}, \phi \in \mathcal{Z} \label{for_drl_scale:lm} \\
    & \textit{vm}^{\phi}_{p} = \sum_{d \in {D}} \textit{vd}_d \sum_{k \in {K}} w_k \sum_{\textit{tr} \in \textit{TR}^\textit{OB}(p)} \sum_{\textit{bl} \in \textit{BL}} \sum_{b \in{B}} \Tilde{u}^{b,d,\textit{bl},\phi}_{\textit{tr},k} \nonumber\\
    & \qquad\qquad\qquad \forall p \in {P}, \phi \in \mathcal{Z} \label{for_drl_scale:vm} \\
    & \underline{\textit{lcg}} \cdot \textit{tw}^{\phi}_{p} \leq \textit{lm}^{\phi}_{p} \leq \overline{\textit{lcg}} \cdot \textit{tw}^{\phi}_{p} \nonumber\\
    & \qquad\qquad\qquad \forall p \in {P}, \phi \in \mathcal{Z} \label{for_drl_scale:LS1} \\
    & \underline{\textit{vcg}} \cdot \textit{tw}^{\phi}_{p} \leq \textit{vm}^{\phi}_{p} \leq \overline{\textit{vcg}} \cdot \textit{tw}^{\phi}_{p} \nonumber\\
    & \qquad\qquad\qquad \forall p \in {P}, \phi \in \mathcal{Z} \label{for_drl_scale:VS1} \\
    & \Tilde{u}^{b,d,\textit{bl},\phi'}_{\textit{tr},k} = \Tilde{u}^{b,d,\textit{bl},\phi}_{\textit{tr},k} \nonumber\\
    & \qquad\qquad\qquad \forall p \in {P}, \textit{tr} \in \textit{TR}^{+}(p), k \in {K}, b \in{B}, d \in {D}, \textit{bl} \in \textit{BL}, \nonumber\\
    & \qquad\qquad\qquad  \phi, \phi' \in \mathcal{Z} \mid q^{\phi}_{[p-1]} = q^{\phi'}_{[p-1]} \label{for_drl_scale:nonanti}
\end{align}
\section{Deep Reinforcement Learning Framework} \label{sec_drl_scale:DRL_model}
In the previous section, an explicit MIP model is provided that grows rapidly with problem size due to time-indexed variables and constraints. As an alternative, we propose an MDP formulation of the MPP, modeling it as a sequential decision process. This structure scales with episode length and the sizes of the state and action spaces, avoiding the combinatorial growth inherent in MIPs. In this section, we extend the DRL framework from \cite{van_twiller_navigating_2025} by refining the MDP, the policy model and feasibility mechanisms. For completeness, we restate the formal and decomposed MDP formulations for the MPP under demand uncertainty, incorporating extensions for blocks and PBS patterns. We then introduce the architecture of AI2STOW, integrating updated self-attention (SA) mechanisms into the dynamic embeddings, before defining the action mask to enforce PBS.

\subsection{Formal Markov Decision Process}
We define an episodic discounted MDP to model the MPP under demand uncertainty. The MDP is represented as \(\mathcal{M} = (S, X, \mathcal{T}, \mathcal{R}, P^{N_P-1}_1, \gamma), \) where $T_\textit{epi} = N_1^{N_P-1}$ is a finite horizon equal to the number of load ports. Let us also define the shapes of the utilization, vessel location, and cargo demand as \( n_u = |B| \times |D| \times |\textit{BL}| \times |K| \times |\textit{TR}| \), \( n_c = |B| \times |D| \times |\textit{BL}| \), and \( n_q = |K| \times |\textit{TR}| \), respectively. 

Before formally introducing the aspects of \( \mathcal{M} \), we first provide some intuition on this sequential decision-making problem. Each step of the MDP corresponds to a port call, where the state captures the vessel’s utilization $u_p$ and the cargo demand at each port $q_p$. An action $x_p$ represents the placement of specific cargo at vessel locations in the current port. Transitions related to vessel utilization involve adding loaded cargo and removing discharged cargo, while transitions related to demand correspond to the realization of demand upon port arrival. The reward at each step is the cargo revenue minus the hatch overstowage and excess crane move costs.

Formally, the current state \( s_p \in S \) is defined as \( s_p = (u_p, q_p, \zeta) \), where vessel utilization \( u_p \in \mathbb{R}_{\geq0}^{n_u} \), and realized demand \( q_p \in \mathbb{R}_{\geq0}^{n_q} \). The environment parameters \( \zeta \) include the expected demand \( \mu \in \mathbb{R}_{\geq 0}^{n_q} \) and its standard deviation \( \sigma \in \mathbb{R}_{>0}^{n_q} \). The demand is characterized by load ports \( i \), discharge ports \( j \), and cargo types \( k \), where \( (i,j,k) \in \textit{TR} \times K \). Further, the TEU per container is given by \( \textit{teu} \in \{1,2\}^{n_q} \), cargo weight by \( w \in \mathbb{R}_{>0}^{n_q} \), and cargo revenue by \( \textit{rev} \in \mathbb{R}_{>0}^{n_q} \). For consistency in dimensionality, \( \textit{teu} \) and \( w \), previously defined with shape \( |K| \) in Section~\ref{sec_drl_scale:problem}, are expanded to shape \( n_q = |K| \times |\textit{TR}| \) in the MDP. The initial state is given by \( s_0 = (u_0, q_0, \zeta) \), where the vessel starts empty with \( u_0 = \mathbf{0}^{n_u} \), and the realized demand \( q_0 \) at the initial port. The parameters \( \zeta \) are randomly initialized at the start of each episode.

The action \( x_p \in X \) represents the placement of containers on the vessel, updating the utilization \( u_p \), where \( x_p \in \mathbb{R}_{>0}^{n_u} \) and is analogous to \(u_{\textit{tr},k}^{b,d,\textit{bl}} \; \forall \textit{tr} \in \textit{TR}^+_p, k \in K, b \in B, d \in D, \textit{bl} \in \textit{BL}\). Actions $ {x}_p $ are subject to a feasible region, defined by polyhedron $ \textit{PH}({s}_p) = \{{x}_p \in \mathbb{R}^{n_u}_{>0} : A({s}_p) {x}_p \leq b({s}_p)\} $. Here, $ A({s}_p) \in \mathbb{R}^{m_u \times {n_u}} $ is the constraint matrix, $ b({s}_p) \in \mathbb{R}^{m_u} $ is the bound vector, and $m_u$ is the number of constraints. As shown in \cite{pacino_fast_2011}, linearly relaxed MPPs outperform integer MPPs, as the subsequent SPP discretizes the solution. Consequently, we adopt real-valued actions instead of integer-valued actions. 

Given a traditional MPP with utilization \( u_p \) \cite{van_twiller_efficient_2024}, the constraints of \( \textit{PH}({s}_p) \) are defined in terms of actions \( x_p \) and pre-loading utilization \( u'_p = u_{\textit{tr},k}^{b,d,\textit{bl}} \; \forall \textit{tr} \in \textit{TR}^\textit{ROB}_p, k \in K, b \in B, d \in D, \textit{bl} \in \textit{BL}\). We also introduce the Hadamard product $\odot$ for element-wise multiplication.
\begin{align}
\quad {x}_p^\top \mathbf{1}_{n_c} & \leq \; q_p  \label{for_drl_scale:fr_demand}\\
\quad  {x}_p \textit{teu} & \leq \; c - u'_p \textit{teu} \label{for_drl_scale:fr_capacity}\\
\text{-}\mathbf{1}^\top \! \big((\textit{lm}\!-\!\underline{\textit{lcg}} w)\!\odot\!{x}_p \big) & \leq \! \text{-}\mathbf{1}^\top\! \big( (\underline{\textit{lcg}} w\!-\!\textit{lm})\!\odot\!u'_p \big) \label{for_drl_scale:fr_lcg_lb} \\
\;\mathbf{1}^\top \! \big( (\textit{lm}\!-\!\overline{\textit{lcg}} w)\!\odot\!{x}_p \big) & \leq \! \; \mathbf{1}^\top\!\big( (\overline{\textit{lcg}} w\!-\!\textit{lm})\!\odot\!u'_p \big) \label{for_drl_scale:fr_lcg_ub} \\
\text{-}\mathbf{1}^\top \! \big((\textit{vm}\!-\!\underline{\textit{vcg}} w)\!\odot\!{x}_p \big) & \leq \! \text{-}\mathbf{1}^\top\!\big( (\underline{\textit{vcg}} w\!-\!\textit{vm})\!\odot\!u'_p \big) \label{for_drl_scale:fr_vcg_lb} \\
\;\mathbf{1}^\top \! \big((\textit{vm}\!-\!\overline{\textit{vcg}} w)\!\odot\!{x}_p \big) & \leq \! \; \mathbf{1}^\top\!  \big( (\overline{\textit{vcg}} w\!-\!\textit{vm})\!\odot\!u'_p \big) \label{for_drl_scale:fr_vcg_ub} \\
x_p & \leq M ( \textit{di}(u'_p) + 1 - \textit{eu}(u'_p))  \label{for_drl_scale:fr_bs}
\end{align}

Constraint \eqref{for_drl_scale:fr_demand} restricts the load \( x_p \) to the available demand \( q_p \), while Constraint \eqref{for_drl_scale:fr_capacity} ensures that \( x_p \) does not exceed the residual TEU capacity. Constraints \eqref{for_drl_scale:fr_lcg_lb}–\eqref{for_drl_scale:fr_lcg_ub} impose limits on the longitudinal center of gravity (LCG), whereas Constraints \eqref{for_drl_scale:fr_vcg_lb}–\eqref{for_drl_scale:fr_vcg_ub} enforce similar bounds on the vertical center of gravity (VCG). These stability constraints are derived from Constraints \eqref{for_drl_scale:TW}–\eqref{for_drl_scale:VS1}, as presented in \cite{van_twiller_navigating_2025}. Finally, Constraint~\eqref{for_drl_scale:fr_bs} enforces PBS patterns on the action variable \( x_p \), where \( \mathit{di}(u'_p) \in \{0,1\}^{n_u} \) indicates of available locations based on used PODs in \(u'_p\), and \( \mathit{eu}(u'_p) \in \{0,1\}^{n_u} \) denotes currently empty locations based on \(u'_p\).

The stochastic transition function \( \mathcal{T}({s}_{p+1} | {s}_p, {x}_p) \in \Delta(S) \) updates the current state. During episodes, the transition comprises multiple components:
\begin{itemize}
    \item Upon port arrival, port demand \( q_p \) is revealed.
    \item Subsequently, onboard cargo is discharged \( u_{p+1} = u_p \odot (1 - \mathbf{e}_p^-) \), where \( \mathbf{e}_p^- \in \{0, 1\}^{n_q} \) is a binary mask indicating the cargo type and transport to nullify in \( u_p \).
    \item Finally, cargo is loaded onboard \( u_{p+1} = u_p + {x}_p\). Action ${x}_p$ is based on the current utilization $u_p$ and revealed demand $q_p$ of port $p$. Future port demand stays unknown.
\end{itemize}

The deterministic reward function is defined in Equation \eqref{for_drl_scale:reward} and computes profit as the difference between revenue and costs. Revenue is given by \( \min(q, x^\top \mathbf{1}_{n_c})^\top \textit{rev} \), where \( \textit{rev} \in \mathbb{R}^{n_q}_{>0} \) is the per-unit revenue, \( x \in \mathbb{R}^{n_q \times n_c}_{\geq 0} \) the cargo allocation, and \( q \in \mathbb{Z}^{n_q}_{>0} \) the demand vector. The pointwise minimum ensures revenue is only collected up to demand. Costs are determined using state-dependent auxiliary variables: hatch overstows \( \textit{ho}({s}_p, p) \in \mathbb{R}^{|B|}_{>0} \) and excess crane moves \( \textit{cm}({s}_p, p) \in \mathbb{R}^{|B|-1}_{>0} \). These costs are weighted by \( \textit{ct}^\textit{ho} \in \mathbb{R}_{>0} \) and \( \textit{ct}^\textit{cm} \in \mathbb{R}_{>0} \), respectively.
\begin{align} \label{for_drl_scale:reward}
\mathcal{R}(s, x, p) = & \; \textit{rev} \min\big(x^\top \mathbf{1}_{n_c}, q \big)\nonumber \\  & - \textit{ct}^\textit{ho} \mathbf{1}^\top \textit{ho}(s, p) + \textit{ct}^\textit{cm} \mathbf{1}^\top \textit{cm}(s, p).
\end{align}

\subsection{Decomposed Markov Decision Process}
The formal MDP defines an action space of size \( |X| \propto |B| \cdot |D| \cdot |\textit{BL}| \cdot |K| \cdot |P| \). An action simultaneously places all container types and transports on the vessel, which obtains an action space with high dimensionality that may impede learning efficiency \cite{kanervisto_action_2020}.  An alternative is to decompose the formal MDP as a sequence of granular decisions indexed by \( (i,j,k) \), where \( (i,j) \in \textit{TR} \) represents transport routes and \( k \in K \) denotes cargo types. Consequently, actions sequentially handle each transport \( (p,j) \) and cargo type \( k \) at port \( p \), followed by sailing to the next port. This decomposition reduces the action space to \( |X| = |B| \cdot |D| \cdot |\textit{BL}| \) while distributing decisions over an extended time horizon \( t \in H = \{0,1,\dots,T_\textit{seq}\} \) with  
\(T_\textit{seq} = |K| \cdot |\textit{TR}|.\) Due to its reduced size, we prefer the decomposed MDP. 

The state $s_t = (u_t, q_t, \zeta)$ depends on time step $t$, including vessel utilization $u_t \in \mathbb{R}^{n_u}_{\geq 0}$ and realized demand $q_t \in \mathbb{R}^{n_q}_{\geq 0}$. The environment parameter \( \zeta \) remains unchanged. Given the time $t$, however, we can extract relevant parameters from \( \zeta \), such as \(\textit{pol}_t, \textit{pol}_t, k_t, \textit{rev}^{(\textit{pol}_t, \textit{pod}_t, k_t)}\). 

At each step \( t \), the action \( x_t \in \mathbb{R}^{n_c}_{\geq 0} \) assigns containers to the utilization \( u_t \). Each action is constrained by the feasible region \( \textit{PH}(s_t) \), defined as  
\(\textit{PH}(s_t) = \{x_t \in \mathbb{R}^{n_c}_{>0} : A(s_t) x_t \leq b(s_t)\}. \) Here, \( A(s_t) \in \mathbb{R}^{m_c \times n_c} \) is the constraint matrix, \( b(s_t) \in \mathbb{R}^{m_c} \) is the bound vector, and \( m_c \) represents the number of constraints. Note that Constraints \eqref{for_drl_scale:fr_demand}-\eqref{for_drl_scale:fr_bs} are reformulated to align with \( \textit{PH}(s_t) \).

The stochastic transition consists of loading at each time step \(t\), where \( x_t \) is added to \( u_t \), while discharging and demand realization occur only upon arrival at a new port, indicated by \( t \in T_{\textit{new port}} \):
\begin{align}
T_{\textit{new port}} &= \Big\{\!t\!\in\!H\! \ \mid \!\exists p\!\in\!P^{N_P-1}_1, t = |K|\Bigl((p-1)(N_P\!-1\!) - \frac{p(p-1)}{2}\!\Bigr)\!\Big\}.
\end{align}

The revenue at step \( t \) is computed as \( \textit{rev}(\textit{pol}_t, \textit{pod}_t, k_t)\mathbf{1}^\top x_t. \) However, costs depend on knowing all loading operations at port $p$, which is aggregated in utilization $u_t$ at the last step of the port \( t \in T_{\textit{leave port}} \):
\begin{align}
T_{\textit{leave port}} &= \Big\{\!t\!\in H\! \ \mid\! \exists p\! \in\! P^{N_P-1}_1, t = |K|\Bigl(p(N_P\!-1\!) - \frac{p(p-1)}{2}\Bigr) - 1 \Big\}.
\end{align}
 As a result, the cost signal is sparse, being evaluated only once per port $p$ rather than at each step.

\subsection{Proposed Architecture} \label{sec_drl_scale:drl_arch}  
\begin{figure}[h!]
    \centering
    \includegraphics[width=\linewidth]{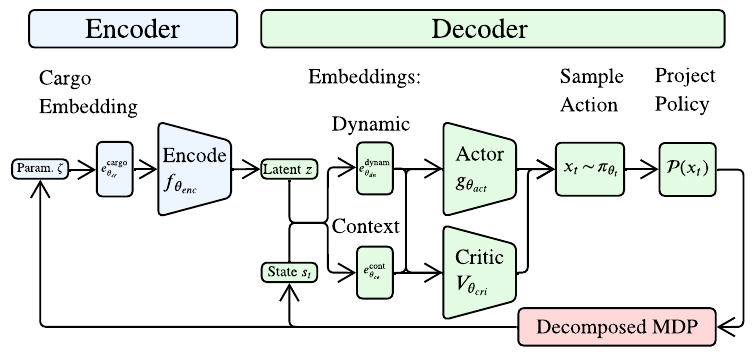}
    \caption[Deep reinforcement learning architecture with feasibility projection]{Deep reinforcement learning architecture with feasibility projection for actor-critic methods \cite{van_twiller_navigating_2025}}
    \label{fig_drl_scale:drl_architecture}
\end{figure}

The AI2STOW architecture comprises three key components: an encoder-decoder model parameterized by \( \theta \), a feasibility projection layer and the DRL implementation. As illustrated in Figure~\ref{fig_drl_scale:drl_architecture}, the encoder-decoder model employs a look-ahead policy \( \pi_{\theta}(x | s_t) \) that considers an estimate of future states, which is conditioned on the state \( s_t \) and parameterized by a mean \( \mu_\theta(s_t) \) and standard deviation \( \sigma_\theta(s_t) \). By training within the decomposed MDP and applying policy projection, the framework iteratively generates actions \( x_t \) to construct feasible solutions.

\subsubsection{Encoder-Decoder Model}
\begin{figure}[h!]
    \centering
    \includegraphics[width=0.75\linewidth]{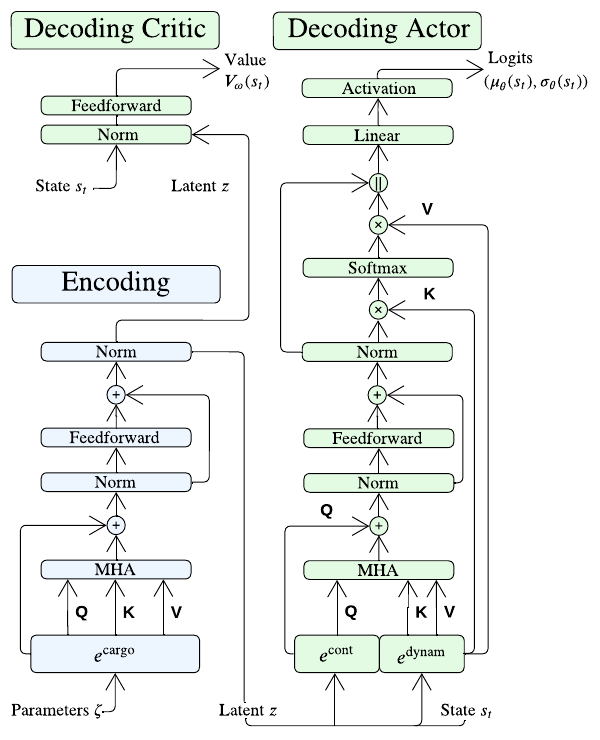}
    \caption[Layers of the encoder and the actor-critic decoder]{Layers of the encoder and the actor-critic decoder \cite{van_twiller_navigating_2025}}
    \label{fig_drl_scale:encoder_decoder_models}
\end{figure}

Figure~\ref{fig_drl_scale:encoder_decoder_models} presents the encoder-decoder model, originally based on the work of \cite{kool_attention_2019}, with architectural adjustments highlighted. Note that $+, \times, \|$ denote summation, multiplication and concatenation operations in Figure~\ref{fig_drl_scale:encoder_decoder_models}.

The core idea behind encoder-decoder models is to convert the input into a compact, low-dimensional representation that allows subsequent models to efficiently extract features and incorporate new information. For instance, the encoder uses embedding $e^\textit{cargo}: 8 \times n_q \rightarrow E$, mapping input $\zeta$ to a representation of shape $E$. The encoder $f: E \rightarrow Z$ maps the embedding output to latent variable $z \in Z$. The context embeddings takes $z$ and utilization $u_t$ as input by $e^\textit{cont}: Z\times n_u \rightarrow E$, whereas the dynamic embedding takes $z$ and $q_t$ as input by $e^\textit{dynam}: Z\times n_q \rightarrow E$. In turn, the actor decoder uses both context and dynamic representations to map $g: 2\times E \rightarrow X$. The embeddings transform different parts of the input data into representations, enabling the model to process specific groups of features and learn meaningful patterns from them. This will be discussed in more detail in the text below.

The cargo embedding \( e^\textit{cargo}(\zeta) \) parameterized by \( \theta_\textit{cr} \) maps cargo-related episode information in $\zeta$ into a feature representation for the encoder. To incorporate positional information, we apply sinusoidal positional encoding \cite{vaswani_attention_2017}, which augments the embeddings with continuous signals that enable the model to distinguish element order and capture position-dependent structure.

An attention encoder $f(e^\textit{cargo}(\zeta))$ parameterized by ${\theta_\textit{enc}}$  maps $e^\textit{cargo}(\zeta)$ to latent variable $z$ using multi-head attention (MHA) to identify relevant features dynamically \cite{vaswani_attention_2017}. 
In MHA, query (Q), key (K), and value (V) are learned projections of the input that enable the model to selectively weigh and aggregate input elements based on their contextual relevance, where attention weights are computed by comparing queries with keys to determine how much focus to place on each value. Then, we use a feed-forward network (FFN) with ReLU activation to introduce non-linearity, layer normalization to stabilize training, residual connections to sum the MHA output with its input to facilitate gradient flow and preserve input information, and dropout to prevent overfitting \cite{goodfellow_deep_2016}.

The context embedding $e^\textit{cont}(u_t, z)$, parameterized by ${\theta_\textit{co}}$, focuses on time $t$, extracting time-specific features from the utilization $u_t$ and latent variable $z$. These features provide the MHA query with a representation of the vessel’s current condition, enabling the policy to make decisions based on the present.

The dynamic embedding \(e^\textit{dynam}(q_t, z)\), parameterized by \(\theta_\textit{dn}\), leverages SA to capture temporal dependencies in demand \(q_t\) over horizon \(H\). By dynamically weighing past and present observations, SA identifies trends and long-term patterns, enhancing decision-making \cite{vaswani_attention_2017}. This embedding combines real demand and latent factors to produce keys and values for an MHA layer. By emphasizing relevant historical patterns while filtering noise, it enables the model to anticipate future conditions and optimize long-term strategies.

As proposed in \cite{kool_attention_2019}, our actor decoder $g(e^\textit{cont}(u_t, z), e^\textit{dynam}(q_t, z))$ is also an attention model, parameterized by $\theta_\textit{act}$. However, our actor takes input from the context and dynamic embedding. An MHA layer with a forward-looking mask bases decisions on steps ${t, t+1, \dots, T_\textit{seq}}$ to anticipate future events. Subsequently, an FFN extracts features, after which a pointer mechanism performs a soft selection of relevant steps using key-value inputs \cite{vinyals_pointer_2015}. A softplus activation outputs positive logits $({\mu}_\theta(s_t), {\sigma}_\theta(s_t))$.

Our critic model \( V_{\theta_\textit{cri}}(s_t, z) \), parameterized by \( \theta_\textit{cri} \), estimates the value of state \( s_t \) and latent variable \( z \) through an FFN outputting \( V_{\theta_\textit{cri}}(s_t, z) \in \mathbb{R} \), providing a low-variance baseline to guide and stabilize policy updates.

The actor logits parameterize a stochastic policy \( \pi_{\theta}(x \mid s_t) = \mathcal{N}(\mu_\theta(s_t), \sigma_\theta(s_t)) \), which enables sampling actions \( x_t \sim \pi_{\theta}(x \mid s_t) \) to construct solutions, while promoting exploration and supporting gradient-based optimization.

\subsubsection{Feasibility Mechanisms}
A variety of feasibility mechanisms can be used to ensure policy $\pi_\theta$ adheres to feasible region $\textit{PH}(s_t)$. A straightforward approach is to include feasibility regularization (FR) in the loss function, as shown in Equations \eqref{for_drl_scale:actor_loss} and \eqref{for_drl_scale:feas_loss}. In previous work, we argued that FR faces significant challenges when applied to dynamic, state-dependent feasible regions \cite{van_twiller_navigating_2025}. The primary difficulty lies in determining the Lagrangian multiplier \(\lambda_f\), as balancing multiple constraints with varying scales requires expensive hyperparameter optimization. Despite these limitations, FR remains a useful baseline for comparison.
\begin{align} 
\mathcal{L}(\theta) = & -\mathcal{L}_{\textit{actor}}(\theta) + \lambda_{{f}} \mathcal{L}_{\textit{feas}}(\theta)  \label{for_drl_scale:actor_loss} \\
\mathcal{L}_{\textit{feas}}(\theta) = &  \mathbb{E}_t \left[ \left( A(s_t) x_\theta(s_t) - b(s_t) \right)_{>0} \right] \label{for_drl_scale:feas_loss}
\end{align}

Projection layers \(\mathcal{P}: X \rightarrow X\) enforce feasibility by mapping policy samples onto feasible sets by non-linear transformations, which distort the original action density and require Jacobian-based log-probability corrections~\cite{bishop_pattern_2006}. However, such corrections can be undefined, leading to biased gradients, particularly when using closed-form (convex) solvers like~\cite{agrawal_differentiable_2019}. Moreover, these solvers are often too slow for practical use during training. We instead favor efficient, differentiable projection layers with well-defined Jacobians, such as the violation projection (VP) layer in~\cite{van_twiller_efficient_2024}.  Regardless of the limitations, solvers remain valuable at inference, where differentiability and speed are less critical, providing a reliable post-processing step to enforce feasibility.

Furthermore, \(\textit{PH}(s_t) \) transforms from a convex to a non-convex polyhedron if PBS is included, as the introduction of binary values creates a discontinuity in the polyhedron. To preserve convexity, we propose omitting Constraint \eqref{for_drl_scale:fr_bs} from \( \textit{PH}(s_t) \) to allow VP to handle the remaining convex constraints, whereas the non-convex constraint can be handled by action masking. Masking is a differentiable operation that preserves the linear structure, allowing constraint enforcement without output distortion \cite{stolz_excluding_2025}. We define $\textit{xm}(s_t) \in \{0,1\}^{n_q}$ as a state-dependent action mask. However, since the policy essentially says not to execute an action for some elements in $x_t$, we want the probability of this element occurring to be 0. Hence, Equations \eqref{for_drl_scale:action_masking} and \eqref{for_drl_scale:logprob_masking} define the action masking and handling of log probabilities.
\begin{align}
    x_t' &= \textit{xm}(s_t) \odot x_t \label{for_drl_scale:action_masking} \\
    \log \pi(a | s) &= \textit{xm}(s_t) \odot \log \pi(a | s) + (1 - \textit{xm}(s_t)) \cdot (-\infty) \label{for_drl_scale:logprob_masking}
\end{align}

Subsequently, we define a state-dependent action mask \(\textit{xm}(s_t)\) to enforce PBS patterns. Algorithm~\ref{alg_drl_scale:PBS_mask} defines the PBS action mask, which ensures the grouping of cargo with POD \( j \) into blocks while maximizing the available space. Recall that blocks include both above-deck and below-deck locations. The algorithm works as follows:
\begin{enumerate}
    \item The algorithm starts with computing POD demand, residual capacity and POD locations indicator based on the state $s_t$. Then, it identifies locations that are either empty or already used by POD \( j \).
    \item The algorithm calculates the amount of demand required to be placed in empty locations, given the residual capacity.
    \item To maintain stability, a random score is assigned to half of the bay-block pairs and reflected symmetrically to the other half. Only empty bay-block pairs are considered.
    \item The scored bay-block pairs are sorted in descending order, and a top-\(k\) subset is selected such that the cumulative capacity of those empty locations is sufficient to meet the remaining demand.
    \item Finally, the mask is constructed by marking the selected top-\(k\) empty locations as valid, while preserving locations already used by POD \( j \).
\end{enumerate}

\begin{algorithm}[h!]
\caption{Action Mask for PBS}
\label{alg_drl_scale:PBS_mask}
\small
\begin{algorithmic}[1]
\Require POD $j \in P$, realized demand $q \in \mathbb{Z}^{n_q}_{\geq 0}$, TEU capacity $c \in \mathbb{Z}^{n_c}_{>0}$, utilization $u \in \mathbb{R}^{n_c \times n_q}_{\geq 0}$
\Ensure Valid location mask $\mathit{xm} \in \{0,1\}^{n_c}$
\State \textbf{Compute Utilization Status:}
\State POD demand: $q_{\textit{pod}} = \sum_{i \in P} \sum_{k \in K} q^{(i,j,k)}$
\State Residual capacity: $c_{\textit{res}} = c - u \cdot {\textit{teu}}$
\State POD indicator: $\textit{il}_{b,d,\textit{bl},j} = \mathbb{I}\left( \sum_{i \in P} \sum_{k \in K} u^{(b,d,\textit{bl},i,j,k)} > 0 \right) \forall b \in B,\ d \in D,\ \textit{bl} \in \textit{BL},\ j \in P$
\State Empty locations: $\textit{el} = \{(b,d,\textit{bl}) \mid \textit{il}_{b,d,\textit{bl},j} = 0, \forall j \in P \}$
\State Used locations with POD $j$: $\textit{ul}_j = \{(b,d,\textit{bl}) \mid \textit{il}_{b,d,\textit{bl},j} = 1 \}$ \\

\State \textbf{Compute Remaining Demand to Fulfill:}
\State $\textit{rq} \gets \max\left(0,\ q_{\textit{pod}} - \sum_{(b,d,\textit{bl}) \in \textit{ul}_j} c_{\textit{res}}^{b,d,\textit{bl}} \right)$
\State $\textit{fq} \gets \min\left( \mathbf{1}^\top c_{\textit{res}},\ \textit{rq} \right)$\\

\State \textbf{Symmetric Random Scoring of Empty Bay Blocks:}
\State Sample random scores: $\textit{sc}_{b,\textit{bl}} \sim \mathcal{U}(0,1)$
\State Reflect symmetrically: $\tilde{\textit{sc}}_{b,\textit{bl}} = \textit{sc}_{\min(b,\ |B|-b+1),\ \textit{bl}}$
\State Bay-block empty mask: $\textit{ms}_{b,\textit{bl}} = \bigwedge_{d \in D} \mathbb{I}\left[(b,d,\textit{bl}) \in \textit{el} \right]$
\State Apply mask: $\tilde{\textit{sc}} \gets \tilde{\textit{sc}} \odot \textit{ms}$
\State Sort indices: $\Pi \gets \text{argsort}(\text{vec}(\tilde{\textit{sc}}),\ \text{descending})$ \\

\State \textbf{Select Top-$k$ Empty Locations Meeting Capacity Requirement:}
\State $c_{\Pi} \gets \text{gather}(c\ \text{by}\ \Pi)$
\State $\hat{c} \gets \text{cumsum}(c_{\Pi})$
\State $k \gets \min \left\{ i\ |\ \hat{c}_i \geq \textit{fq} \right\}$
\State $\tilde{\mathit{xm}}[i] \gets \mathbb{I}[i \in \{\Pi_1, \dots, \Pi_k\}]$ \\

\State \textbf{Merge Mask with Used Locations for POD:}
\State $\mathit{xm} \gets \tilde{\mathit{xm}} \land \textit{el} \lor \textit{ul}_j$ \\
\Return $\mathit{xm}$
\end{algorithmic}
\end{algorithm}

\subsubsection{Deep Reinforcement Learning Implementation}
Several DRL algorithms can be employed to train a policy in conjunction with a feasibility projection layer. To this end, we provide a general pseudocode outlining the core steps of the training process in Algorithm~\ref{alg_drl_scale:training}. This pseudocode abstracts the common structure of DRL training loops. The training involves iterative interactions between the agent and the environment. Within the interactions, the policy produces an action based on the current state, which is first element-wise multiplied with state-dependent mask $\textit{xm}(s_t)$, and then passed through a projection layer for convex constraints. The resulting projected actions are used to collect transitions, which are stored in an experience buffer. This buffer is subsequently used to update value estimates, compute the loss function that guides policy improvement, and update the model parameters $\theta$. The pseudocode is algorithm-agnostic and can be instantiated as either on-policy (e.g., PPO) or off-policy (e.g., SAC), depending on how the experience buffer $\mathcal{D}$ is populated and sampled.

\begin{algorithm}[h!]
\caption{General DRL Training with Feasibility Projection}
\label{alg_drl_scale:training}
\small
\begin{algorithmic}[1]
\Require MDP $\mathcal{M}$, initial policy $\pi_\theta$, projection layer $\mathcal{P}(\cdot)$, mask layer $\textit{xm}(\cdot)$, buffer $\mathcal{D}$
\Ensure Trained policy $\pi_\theta$

\For{each training iteration}
    \State Interact with environment $\mathcal{M}$ using masked and projected policy $\mathcal{P}(\textit{xm}(s_t) \odot \pi_\theta(s_t))$
    \State Store collected transitions in experience buffer $\mathcal{D}$
    \State Update performance metrics (e.g., $G_t,\ V(s_t),\ Q(s_t,x_t)$) from buffer $\mathcal{D}$
    \State Compute loss $\mathcal{L}(\theta)$ based on buffer $\mathcal{D}$
    \State Update parameters: $\theta \gets \theta - \alpha \nabla_\theta \mathcal{L}(\theta)$
\EndFor

\State \Return $\pi_\theta$
\end{algorithmic}
\end{algorithm}

We also present a general pseudocode for inference in Algorithm~\ref{alg_drl_scale:inference}. During inference, the trained policy generates actions based on the current state. A raw action is sampled from the learned distribution and then passed through the projection layer to enforce feasibility. The projected action is used to compute the reward and transition the environment to the next state. This process is repeated for a fixed time horizon, and the trajectory of visited states and raw actions is recorded for analysis. 

\begin{algorithm}[h!]
\caption{General DRL Inference with Feasibility Projection}
\label{alg_drl_scale:inference}
\small
\begin{algorithmic}[1]
\Require MDP $\mathcal{M}$, trained stochastic policy $\pi_\theta$, projection layer $\mathcal{P}(\cdot)$, mask layer $\textit{xm}(\cdot)$
\Ensure Trajectory $\{(s_t, \tilde{x}_t)\}_{t=0}^{T_\textit{seq}}$, episodic reward $R_{\textit{seq}}$

\State Sample initial state $s_0$ from environment $\mathcal{M}$
\State Initialize episodic reward $R_{\textit{seq}} \gets 0$

\For{$t = 0$ to $T_\textit{seq}$}
    \State Sample raw action $\tilde{x}_t \sim \pi_\theta(s_t)$
    \State Apply masking: $x'_t \gets \textit{xm}(s_t) \odot \tilde{x}_t $
    \State Apply convex projection: $x_t \gets \mathcal{P}(x'_t)$
    \State Observe reward $r_t = \mathcal{R}(s_t, x_t)$
    \State Observe next state $s_{t+1} \sim \mathcal{T}(\cdot \mid s_t, x_t)$
    \State Accumulate reward: $R_{\textit{seq}} \gets R_{\textit{seq}} + r_t$
    \State Update state: $s_t \gets s_{t+1}$
\EndFor

\State \Return $\{(s_t, \tilde{x}_t)\}_{t=0}^{T_\textit{seq}}$, $R_{\textit{seq}}$
\end{algorithmic}
\end{algorithm}

\section{Experimental Results} \label{sec_drl_scale:results}
In this section, we examine the performance of AI2STOW on the MPP under demand uncertainty based on a realistic problem to assess its scalability. We start by discussing the experimental setup, after which we analyze the performance of the AI2STOW policy. Subsequently, an ablation study is performed, after which managerial insights are provided.

\subsection{Experimental Setup}
The largest container vessels range between 14,500 and 24,000 TEU \cite{van_twiller_literature_2024}. To adequately reflect such vessels, we assume the use of a vessel of 20,0000 TEU with 20 bays and 3 hatch covers in each bay.  Typically, stowage planners use a horizon of 5 ports to mitigate the complexity of planning entire voyages. This provides a focused and manageable scope for optimization. These and additional MPP parameters are provided in \Cref{app_drl_scale:mpp_parameters}

In \Cref{app_drl_scale:instance_gen}, we define the instance generator used for the training and evaluation of AI2STOW, including a description of the container and TEU demand per port. An upper bound \(\mathit{ub} \in \mathbb{R}^{n_q}\) is computed to align the expected demand with the vessel's maximum capacity, such that \(\mathbf{1}^\top \mathit{ub} \propto \mathbf{1}^\top c\). During training, \(\mathit{ub}\) is perturbed by a small random factor to encourage instance diversity. Demand instances are then sampled from a uniform distribution \(\mathcal{U}(0, \mathit{ub}^{(i,j,k)})\), for all \((i,j) \in \mathit{TR},\; k \in K\). To test generalization, the voyage length \(N_P\) is varied during evaluation.

Several feasibility mechanisms are used in these experiments to address the constraints in $\textit{PH}(s_t)$. The baseline approach is training DRL with FR (DRL-FR), whereas AI2STOW trains with FR, the general VP layer for convex constraints \cite{van_twiller_navigating_2025}, and the action mask to enforce PBS. Accordingly, the training projection layer, $\mathcal{P}_\textit{train}$, can be denoted as \text{PBS/VP} or similar variants. During inference, a convex programming (CP) layer can be used \cite{agrawal_differentiable_2019}, possibly combined with PBS and policy clamping on TEU capacity (PC) \cite{van_twiller_navigating_2025}. The resulting inference-time projection layer, $\mathcal{P}_\textit{test}$, can thus be expressed as \text{PBS/CP/PC} or a corresponding variant. Implementation details on projection layers are provided in \Cref{app_drl_scale:DRL}.  

We compare AI2STOW against two stochastic programming approaches: a stochastic mixed integer program without anticipation (SMIP-NA), as described in \Cref{sec_drl_scale:problem}, and a stochastic MIP with perfect information (SMIP-PI), which relaxes the non-anticipativity constraint. The SMIP-NA serves as a baseline for AI2STOW, while SMIP-PI represents an expected upper bound. Since AI2STOW operates with continuous actions, we construct the SMIP models using a linear relaxation of the decision variables $\tilde{u}, \tilde{\textit{ho}}$ and $\tilde{\textit{cm}}$.

Training occurs offline on simulated instances provided by the instance generator. Inference (or testing) occurs online on seeded simulated instances: AI2STOW models perform multiple inference rollouts to greedily select the best-performing solution. GPU-based experiments use an NVIDIA RTX A6000, and CPU-based runs use an AMD EPYC 9454 48-Core Processor. Additionally, CPU experiments are given a 1-hour runtime limit, after which the best solution is returned. Implementation details on the DRL algorithms and hyperparameters are also provided in \Cref{app_drl_scale:DRL}.  

\subsection{Policy Performance}
\Cref{tab_drl_scale:comparison} presents a comparative evaluation of AI2STOW against baseline methods. Across all test instances, each version of AI2STOW consistently outperforms SMIP-NA in both objective value, achieving a 20 to 25 \% improvement, and computation time, with reductions of approximately 20 times and 150 times for the CP and VP variants, respectively. With regard to feasibility, only PBS/VP does not fully guarantee feasibility. However, this limitation can be resolved either by fine tuning the parameters in PBS/VP\textsuperscript{\ding{72}} or by applying PBS/VP/PC to clip capacity violations. While these results demonstrate significant performance gains, the SMIP-PI model suggests that further improvements may be attainable if more accurate or complete information were available. Nonetheless, AI2STOW stands out for its efficiency and ability to make high-quality decisions under uncertainty. It also significantly outperforms DRL-FR, achieving a 70 \% higher objective value while maintaining consistent feasibility, whereas DRL-FR failed to produce any feasible solutions. This highlights the practical advantage of projection layers over feasibility regularization.

Among the AI2STOW variants, PBS/CP demonstrates robustness by ensuring feasibility across all instance sizes. However, this comes at the cost of computational efficiency, with a 7 times increase in runtime compared to PBS/VP variants. In contrast, the feasibility performance of PBS/VP and PBS/VP\textsuperscript{\ding{72}} degrades when generalizing to larger instances. This degradation can be attributed to the hyperparameter sensitivity of the VP layer. Larger instances change the episode length and the magnitude of actions, which requires fine-tuning of hyperparameters. To address this, clipping on TEU capacity in PBS/VP/PC achieves feasible solutions with only a slight reduction in objective value compared to PBS/CP. This shows that VP, unlike CP, requires targeted fine-tuning across voyage lengths to be generally effective, while adding PC offers a simple and effective way to increase robustness.

In \Cref{tab_drl_scale:comparison}, there is a lack of results for SMIP models if $N_P\in\{5,6\}$. The exponential growth of the scenario tree with increasing $|\mathcal{Z}| \propto N_P$ significantly impacts computational cost and memory usage. While solving SMIP-NA and SMIP-PI for $N_P = 4$ is possible with a computational time of around 30 minutes, longer voyages become impractical on our hardware due to excessive memory demands, as shown in \Cref{app_drl_scale:add_experiments}. Even if we could hold the scenario tree in memory, we would also face intractable computational times for SMIP-NA. This underscores the scalability of DRL models, which can handle uncertainty implicitly. In contrast, all AI2STOW models across all voyage lengths remain well within the practical tractability threshold of 10 minutes for decision support in stowage planning \cite{pacino_fast_2011}.

\begin{landscape}
\begin{table}[h!]
\centering
\footnotesize
\caption[Experimental results comparing the AI2STOW with baselines]{Experimental results comparing the AI2STOW framework with various inference-time projection layers ($\mathcal{P}_\textit{test}$) against baseline and upper-bound models. The DRL models are trained on voyages with 4 ports ($N_P = 4$), while both DRL and SMIP models are evaluated on $N = 30$ test instances. SMIP-NA (non-anticipative) serves as a baseline, while SMIP-PI (assumes perfect information) represents an expected upper bound; both SMIPs use $S_\textit{ST} = 20$ and are solved using CPLEX. Reported metrics include average objective value in profit (O.), inference time in seconds (T.), and the percentage of feasible instances (F.). Generalization performance is assessed on longer voyages with $N_P \in {5, 6}$, beyond the training distribution.} 
\small
\begin{tabular}{llrrrrrrrrr}
\toprule
\multicolumn{2}{c}{\textbf{Methods}} & \multicolumn{3}{c}{\textbf{Testing ($\boldsymbol{N_P = 4}$)}} & \multicolumn{3}{c}{\textbf{Gen. ($\boldsymbol{N_P=5}$)}} & \multicolumn{3}{c}{\textbf{Gen. ($\boldsymbol{N_P=6}$)}} \\
\cmidrule(r){1-2} \cmidrule(r){3-5} \cmidrule(r){6-8} \cmidrule(r){9-11}
\textbf{Model} & $\boldsymbol{\mathcal{P}_\textit{test}(\cdot)}$ & \textbf{O. (\$)} & \textbf{T. (s)} & \textbf{F. (\%)} & \textbf{O. (\$)} & \textbf{T. (s)} & \textbf{F. (\%)} & \textbf{O. (\$)} & \textbf{T. (s)} & \textbf{F. (\%)} \\ 
\midrule
AI2STOW & PBS/CP & 25637.52 & 76.37 & 100.00 & 34840.63 & 105.19 & 100.00 & 43496.30 & 182.63 & 100.00  \\ 
    & PBS/VP & 25352.06 & 10.59 & 86.67 &  34326.75 & 23.90 & 0.00 & 42515.64 & 38.79 & 0.00  \\
    & PBS/VP\textsuperscript{\ding{72}} & 25404.68 & 11.54 & 100.00 &  34326.75 & 23.51 & 0.00 & 42515.64 & 38.34 & 0.00  \\ 
    & PBS/VP/PC & 25346.14 & 10.33 & 100.00 & 34052.90 & 22.83 & 100.00 & 41979.11 & 38.09 & 100.00  \\
\midrule
DRL-FR & - & 14959.52 & 3.19 & 0.00 &  19417.46 & 5.36 & 0.00 &  23471.99 & 7.96 & 0.00 \\ 
SMIP-NA & - & 20647.22 & 1576.77 & 100.00 & - & - & - &  - & - & - \\ 
SMIP-PI\textsuperscript{*} & - & 33787.63 & 40.91 & 100.00 & - & - & - &  - & - & - \\
\bottomrule
\end{tabular}
\label{tab_drl_scale:comparison} 
\end{table}

\end{landscape}

\subsection{Ablation Study}
\Cref{tab_drl_scale:ablation} presents the results of an ablation study, which replaces or removes model components of the model to assess their impact on performance. Replacing SA with FF layers in the dynamic embedding, and substituting the AM policy with an MLP, led to only marginal performance degradation when using PBS/CP as postprocessing step in  $\mathcal{P}_\textit{test}$. However, under PBS/VP projection, as used during training, replacing SA with FF layers resulted in more feasibility violations. Interestingly, replacing the AM with an MLP under PBS/VP improved feasibility performance. This suggests that the problem instances are sufficiently regular for simpler architectures to learn (e.g., MLPs). However, AMs are likely to outperform MLPs in less structured and more challenging instances, where long-range dependencies and complex interactions are present \cite{vaswani_attention_2017}.

Regarding the projection layers, training with VP works well when using PBS/CP during inference, even improving the objective at the cost of longer runtimes. However, when applying PBS/VP\textsuperscript{\ding{72}} or VP\textsuperscript{\ding{72}} during inference, feasibility drops significantly, as PBS patterns are not enforced. This indicates that PBS/CP can effectively resolve feasibility issues found during training. Still, we prefer solutions to remain closer to the feasible region and rely less on post-processing through $\mathcal{P}_\textit{test}(\cdot)$. Notably, training solely on PBS leads to significant reductions in the objective value and depends on post-processing for feasibility.

\begin{table}[h!]
\centering
\captionsetup{font=footnotesize}
\caption[Ablation study of AI2STOW components]{{Ablation study results for variations of AI2STOW components on $N=30$ instances. The table specifies the DRL algorithm used for training (Alg.), the policy model as an attention model (AM) or a multilayer perceptron (MLP), combined with a dynamic embedding layer using either self-attention (SA) or a feedforward layer (FF), along with the training projection $\mathcal{P}_\textit{train}$ and testing projection $\mathcal{P}_\textit{test}$. Reported metrics include the average objective value in profit (O.), inference time in seconds (T.), and the percentage of feasible instances (F.).}}
\label{tab_drl_scale:ablation}
\small
\begin{tabular}{llllrrr}
\toprule
\multicolumn{4}{c}{\textbf{Methods}} & \multicolumn{3}{c}{\textbf{Testing ($\boldsymbol{N_P=4}$)}} \\
\cmidrule(r){1-4} \cmidrule(r){5-7}
\textbf{Alg.} & \textbf{Model} & $\boldsymbol{\mathcal{P}_\textit{train}(\cdot)}$ & $\boldsymbol{\mathcal{P}_\textit{test}(\cdot)}$  & \textbf{O. (\$)} & \textbf{T. (s)} & \textbf{F. (\%)} \\
\midrule
SAC & SA-AM & PBS/VP & PBS/CP & 25637.52 & 76.37 & 100.00  \\
SAC & SA-AM & PBS/VP & PBS/VP & 25352.06 & 10.59 & 86.67  \\
SAC & FF-AM & PBS/VP & PBS/CP & 25558.75 & 79.86 & 100.00 \\
SAC & FF-AM & PBS/VP & PBS/VP & 25224.73 & 10.51 & 76.67  \\
SAC & FF-MLP & PBS/VP & PBS/CP & 25463.42 & 71.74 & 100.00 \\
SAC & FF-MLP & PBS/VP & PBS/VP & 25222.44 & 9.89 & 96.67 \\
\midrule
SAC & SA-AM & VP & PBS/CP  & 26722.24 & 93.51 & 100.00 \\
SAC & SA-AM & VP & PBS/VP\textsuperscript{\ding{72}}  & 25755.81 & 16.64 & 0.00 \\
SAC & SA-AM & VP & VP\textsuperscript{\ding{72}}  & 30070.51 & 6.36 & 0.00  \\
SAC & SA-AM & PBS & PBS/CP & 10501.61 & 9.76 & 100.00 \\
SAC & SA-AM & PBS & PBS/VP\textsuperscript{\ding{72}}  & 10501.57 &  4.10 & 100.00 \\
SAC & SA-AM & PBS & PBS & 10367.80 & 3.19 & 0.00 \\
\bottomrule
\end{tabular}
\end{table}

\subsection{Managerial Insights}
\begin{figure}[htbp]
\centering

\begin{subfigure}[t]{0.85\textwidth}
\centering
\begin{tikzpicture}
    \begin{axis}[
        width=\textwidth,
        height=0.45\textwidth,
        xlabel={Number of Scenarios},
        ylabel={Profit (\$)},
        legend style={font=\scriptsize, at={(1.05,0.5)}, anchor=west, legend columns=1, align=left},
        tick label style={font=\scriptsize},
        label style={font=\scriptsize},
        grid=both,
        minor tick num=1,
        major grid style={line width=.5pt,draw=gray!50},
        minor grid style={line width=.5pt,draw=gray!30},
        ymajorgrids=true,
        yminorgrids=true,
        xtick={0,4,8,...,64},
        xmin=2, xmax=26,
        ymin=18000, ymax=35500,
        ytick={20000,22000,..., 40000},
        axis y discontinuity=parallel,
        set layers,
    ]
    \addplot[color=darkred, mark=*, mark size=1, line width=.8pt] coordinates {
        (4,24276) (8,22112) (12,21308) (16, 20945) (20,20647) 
    };
    \addlegendentry{SMIP-NA};
    \addplot[color=darkblue, mark=square*, mark size=1, line width=.8pt] coordinates {
        (4,33989) (8,33911) (12,33837) (16,33857) (20,33788) (24,33764)
    };
    \addlegendentry{SMIP-PI};
    \addplot[dashed, color=darkyellow, line width=.8pt] coordinates {
        (0,25637.52) (40,25637.52)
    };
    \addlegendentry{AI2STOW/PBS/CP};

    \addplot[name path=lower_na, draw=none] coordinates {
        (4,23413) (8,21269) (12,20446) (16,20103) (20,19809) 
    };
    \addplot[name path=upper_na, draw=none] coordinates {
        (4,25140) (8,22955) (12,22169) (16,21788) (20,21485) 
    };
    \addplot[darkred!20, opacity=0.65] fill between[of=upper_na and lower_na];

    \addplot[name path=lower_pi, draw=none] coordinates {
        (4,33089) (8,33065) (12,32997) (16,33021) (20,32967) (24,32934) 
    };
    \addplot[name path=upper_pi, draw=none] coordinates {
        (4,34885) (8,34756) (12,34678) (16,34693) (20,34608) (24,34593)
    };
    \addplot[darkblue!20, opacity=0.65] fill between[of=upper_pi and lower_pi];

    \addplot[name path=lower_na, draw=none] coordinates {
        (0,25048.61) (30,25048.61)
    };
    \addplot[name path=upper_na, draw=none] coordinates {
        (0,26226.44) (30,26226.44)
    };
    \addplot[darkyellow!20, opacity=0.70] fill between[of=upper_na and lower_na];
    \end{axis}
\end{tikzpicture}
\caption{Profit with 95\% CI across scenario sizes.}
\label{fig_drl_scale:objective_values}
\end{subfigure}

\vspace{1em}

\begin{subfigure}[t]{0.85\textwidth}
\centering
\begin{tikzpicture}
\begin{axis}[
    width=\textwidth,
    height=0.45\textwidth,
    xlabel={Number of Scenarios and DRL Model},
    ylabel={Computation Time (s, log scale)},
    ymode=log,
    log basis y=10,
    ymin=1, ymax=10000,
    legend style={font=\scriptsize, at={(1.05,0.5)}, anchor=west},
    tick label style={font=\scriptsize},
    label style={font=\scriptsize},
    grid=both,
    minor tick num=1,
    major grid style={line width=.5pt,draw=gray!50},
    minor grid style={line width=.5pt,draw=gray!30},
    bar width=10pt,
    enlarge x limits=0.15,
    symbolic x coords={4,8,12,16,20,24,AI2STOW},
    xtick={4,8,12,16,20,24,AI2STOW},
    xticklabel style={text height=1.5ex},
    nodes near coords align={vertical},
]

\addplot[
    color=darkred,
    mark=*,
    line width=1pt
] coordinates {
    (4, 1.48) (8, 8.00) (12, 49.52) (16, 309.41) (20, 1576.77) (24, 3654.38)
};
\addlegendentry{SMIP-NA}

\addplot[
    color=darkblue,
    mark=square*,
    line width=1pt
] coordinates {
    (4,1.61) (8, 5.51) (12, 12.07) (16, 24.39) (20, 40.91) (24, 63)
};
\addlegendentry{SMIP-PI}

\addplot+[
    mark=triangle*,
    ybar,
    style={fill=custompurple, draw=black},
    mark options={fill=custompurple}
] coordinates {
    (AI2STOW, 327.50)
};
\addlegendentry{AI2STOW Inference}

\addplot+[
    mark=diamond*,
    ybar,
    style={fill=darkyellow, draw=black},
    mark options={fill=darkyellow}
] coordinates {
    (AI2STOW, 251.13)
};

\addlegendentry{AI2STOW Training}

\addplot[
    color=darkgreen,
    line width=1pt
] coordinates {
    (4, 600) (AI2STOW, 600)
};

\addlegendentry{Stowage Timelimit}

\end{axis}
\end{tikzpicture}
\caption[Average computational time across SMIP scenario sizes and AI2STOW]{Average computational time across scenario sizes for SMIP models and AI2STOW. For AI2STOW, the reported time includes both training and inference averaged over 30 instances. Additional instances would further smooth the training time estimate.}
\label{fig_drl_scale:computational_times}
\end{subfigure}

\vspace{1em}

\begin{subfigure}[t]{0.85\textwidth}
\centering
\begin{tikzpicture}
\begin{axis}[
    ybar,
    bar width=8pt,
    width=\textwidth,
    height=0.5\textwidth,
    enlargelimits=0.15,
    xlabel={\textit{UR} Parameter Value},
    ylabel={Profit (\$)},
    xtick={0.7, 0.9, 1.1, 1.3, 1.5},
    legend style={at={(1.02,0.5)}, anchor=west, font=\scriptsize},
    tick label style={font=\scriptsize},
    label style={font=\scriptsize},
    ymin=0,
    ymajorgrids=true,
    bar shift auto,
    grid style={dashed,gray!30},
    error bars/y dir=both,
    error bars/y explicit,
    legend cell align=left,
    set layers 
]

\addplot+[
    style={fill=darkred!60, draw=black},
    error bars/error bar style={black}] plot coordinates {
    (0.7,18903) +- (631,631)
    (0.9,22908) +- (540,540)
    (1.1,25638) +- (590,590)
    (1.3,26919) +- (545,545)
    (1.5,27215) +- (646,646)
};
\addlegendentry{SA-AM}

\addplot+[
    style={fill=darkblue!60, draw=black},
    error bars/error bar style={black}] plot coordinates {
    (0.7,18043) +- (632,632)
    (0.9,22549) +- (544,544)
    (1.1,26160) +- (602,602)
    (1.3,26805) +- (547,547)
    (1.5,27072) +- (633,633)
};
\addlegendentry{FF-AM}

\addplot+[
    style={fill=darkyellow!60, draw=black},
    error bars/error bar style={black}] plot coordinates {
    (0.7,18557) +- (615,615)
    (0.9,22632) +- (508,508)
    (1.1,25463) +- (544,544)
    (1.3,26761) +- (538,538)
    (1.5,26976) +- (589,589)
};
\addlegendentry{FF-MLP}

\end{axis}
\end{tikzpicture}
\caption{Profit with 95\% CI under distributional shift.}
\label{fig_drl_scale:distribution_shift}
\end{subfigure}

\caption[Sensitivity analysis across scenario size and distribution shift]{Sensitivity analysis across scenario size and distributional shift}
\label{fig_drl_scale:comparison_vertical}
\end{figure}
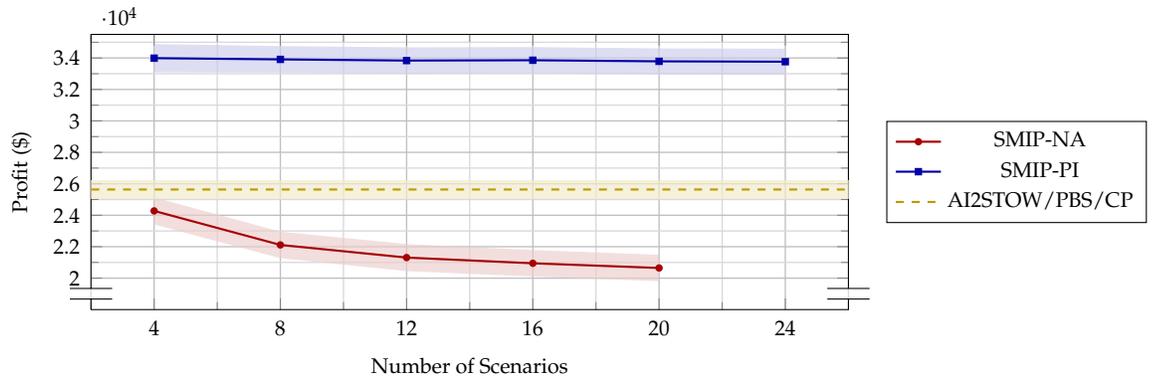
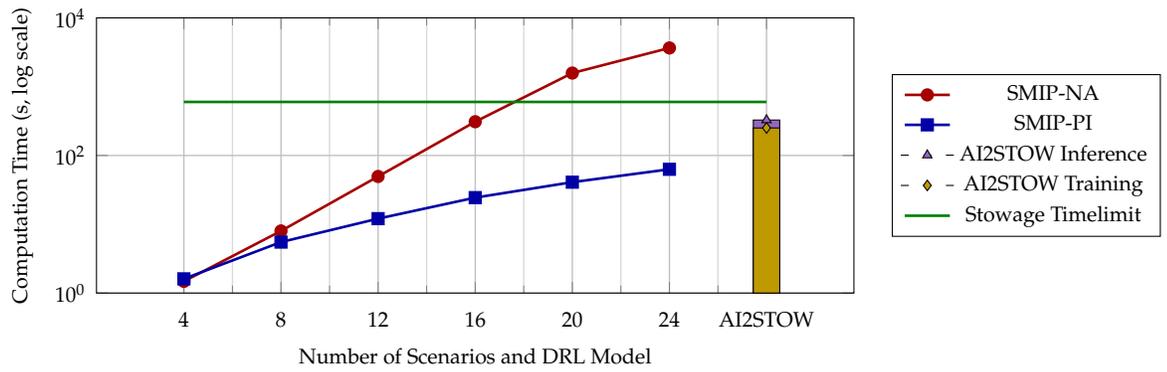
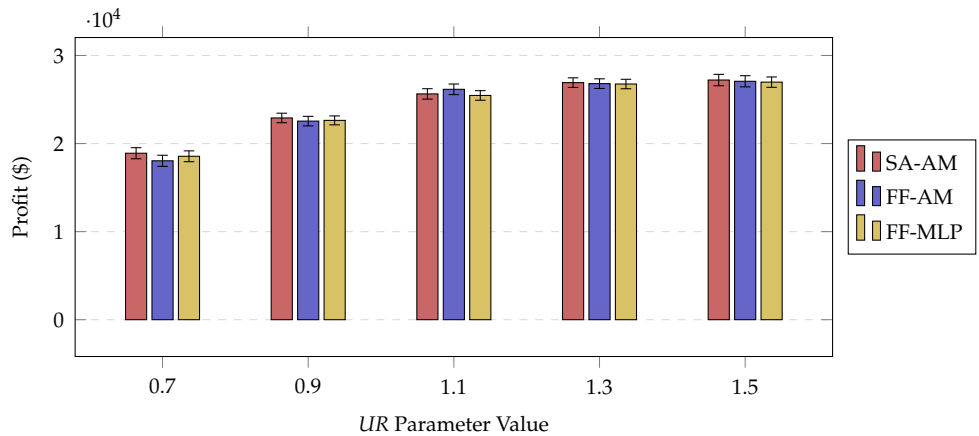

In \Cref{fig_drl_scale:objective_values}, we assess the value of information by comparing profit under non-anticipation (SMIP-NA), imperfect information (AI2STOW), and perfect information (SMIP-PI). Access to additional information leads to substantial improvements in profit relative to SMIP-NA, approximately 25\% for AI2STOW and 60\% for SMIP-PI. While AI2STOW performs significantly better than the non-anticipative baseline, a notable gap remains to the perfect-information upper bound, largely due to the high uncertainty in uniform instances, which are inherently difficult to predict. As the number of scenarios increases, the profit achieved by SMIP-NA gradually decreases and begins to stabilize, indicating diminishing returns from additional information. This trend is less evident in the perfect-information setting. However, convergence is not fully observed within the experimental range: at $S_\textit{ST} = 24$, SMIP-NA returned only empty solutions within the 1-hour time limit, hence, the point is excluded from \Cref{fig_drl_scale:objective_values}. This suggests that beyond a certain scenario size, the scenario tree becomes computationally intractable, regardless of the tractability limit of practical stowage planning.

\Cref{fig_drl_scale:computational_times} illustrates the computational cost of SMIP-NA and SMIP-PI alongside the training and inference time of AI2STOW. Expanding the scenario tree is computationally expensive, with SMIP-NA exhibiting an exponential increase in runtime. The only deviation from this trend occurs at 24 scenarios, where SMIP-NA exceeds the 1-hour time limit. SMIP-PI also shows a steady increase in computational time as scenario size grows. Additionally, the estimated runtime for solving a single instance with AI2STOW is lower than that of SMIP-NA at 20 and 24 scenarios. This suggests that, even with just 30 test instances, the one-time offline training cost of AI2STOW is already justified. As the number of instances increases, this training time is further amortized, improving the cost-effectiveness in larger deployments. If the training budget is sufficient, then the resulting inference time is significantly lower compared to the runtime of SMIP models.

In \Cref{fig_drl_scale:distribution_shift}, the performance of AI2STOW with SA-AM, FF-AM, and FF-MLP is evaluated under varying levels of distributional shift. Note that AI2STOW is trained with $\textit{UR} = 1.1$. The results indicate that all three models perform similarly across different distribution shifts with consistent variability, while SA-AM consistently achieves slightly higher profits in shifting scenarios. This suggests that reducing model complexity does not compromise performance. A closer examination of the distribution shift reveals an asymmetric response. When the distribution is shifted toward lower $\textit{UR}$ values, performance decreases in a manner that scales approximately with the reduced demand, which can be expected as these scenarios fall within the support of the original uniform demand distribution. Conversely, increasing $\textit{UR}$ does not yield improved performance. This indicates that none of the models effectively handle substantial shifts beyond the original distribution, highlighting their limited extrapolation capability under significant demand increases. 

From a practical perspective, we show that including stochastic cargo demand results in a more realistic and complex optimization problem. Particularly challenging to solve using explicit scenario tree formulations, which quickly become intractable. We demonstrate that AI2STOW efficiently generates solutions that implicitly handle uncertainty and consistently identifies master plans with higher expected profit. Even when accounting for training time, our approach is computationally superior to SMIP models. To demonstrate generalization, AI2STOW adapts well to previously unseen and larger voyage lengths. However, under significant distribution shifts, its performance can stagnate. Moreover, its success depends heavily on the quality of training data; when this data is representative of realistic instances, the method performs effectively in practice. This underscores the importance of realistic demand simulators to support the training of DRL-based methods.

\clearpage
\section{Conclusion} \label{sec_drl_scale:conclusion}
This article introduces AI2STOW, an end-to-end deep reinforcement learning model designed to solve the MPP under demand uncertainty, specifically tailored for realistic vessel sizes and operational planning horizons. AI2STOW builds on previous work \cite{van_twiller_navigating_2025} by extending the MDP to include paired block stowage patterns alongside existing global objectives, such as revenue maximization and the minimization of hatch overstowage and excess crane moves. It also respects vessel capacity and ensures both longitudinal and vertical stability.

Experimental results show that AI2STOW generates feasible and adaptive stowage plans for simulated instances on a realistic-sized vessel and operational voyage lengths. AI2STOW also outperforms baseline methods from both deep reinforcement learning and stochastic programming. These findings provide strong evidence that DRL is a promising approach for scalable stowage planning.

Looking ahead, we aim to improve the MDP's representativeness by incorporating local objectives and constraints. We also plan to integrate this framework with the SPP for more efficient end-to-end stowage plan construction. Finally, we encourage exploration into hybrid approaches that combine ML with CO frameworks as alternatives to purely end-to-end learning.

\clearpage
\section*{Acknowledgements}
This work is partially funded by the Danish Maritime Fund (grant nr. 2021-069).

\section*{Declaration of generative AI and AI-assisted technologies in the writing process}
During the preparation of this work, the author(s) used ChatGPT to improve writing. After using this tool/service, the author(s) reviewed and edited the content as needed and take(s) full responsibility for the content of the publication.

\def\urlprefix{}
\def\url#1{}

\bibliographystyle{formatting/elsarticle-harv}
\bibliography{references/references}






\clearpage
\appendix
\section{MPP Parameters} \label{app_drl_scale:mpp_parameters}

The parameters of the MPP are shown in Table \ref{tab_drl_scale:env_params}. 

\begin{table}[h!]
    \centering
    \small
    \caption{MPP parameters} \label{tab_drl_scale:env_params}
    \begin{tabular}{lll}
        \toprule
        \textbf{Parameters} & \textbf{Symbol} & \textbf{Values} \\
        \midrule
        Voyage lengths & $N_P$ & \{4,5,6\} \\
        Number of bays & $N_B$ & 20 \\
        Cardinality deck set & $|D|$ & 2 \\
        Cardinality block set & $|\textit{BL}|$ & 2 \\
        Cardinality cargo set & $|K|$ & 12 \\
        Cardinality transport set & $|\textit{TR}|$ & 6 \\
        Vessel TEU & $\mathbf{1}^\top c$ & 20,000 \\
        Long term contract reduction & $\textit{LR}$ & 0.3 \\
        Utilization rate demand & $\textit{UR}$ & 1.1 \\ 
        lcg bounds & $(\underline{\textit{lcg}},\overline{\textit{lcg}})$ & (0.85,1.05) \\
        vcg bounds & $(\underline{\textit{vcg}},\overline{\textit{vcg}})$ & (0.95,1.15) \\
        Crane moves allowance & $\delta^\textit{cm}$ & 0.25 \\
        Overstowage costs & $\textit{ct}^\textit{ho}$ & 0.33 \\
        Crane move costs & $\textit{ct}^\textit{ho}$ & 0.5 \\
        \bottomrule
    \end{tabular}
\end{table}
\section{Instance Generator} \label{app_drl_scale:instance_gen}
This Appendix introduces the instance generator used to sample demand, after which a descriptive analysis visualizes the demand distribution.

Algorithm \ref{alg_drl_scale:uniform_generator} generates transport matrices of cargo demand based on a perturbed uniform distribution, designed to simulate stowage planning problem instances. Given vessel capacity \( c \) and perturbation factor \( \rho \), an upper bound matrix \( \mathit{ub} \) is first computed proportionally to the total vessel capacity. To introduce controlled randomness, a perturbation is applied element-wise using samples \( U_{i,j,k} \sim \mathcal{U}(0,1) \), scaling each upper bound within the interval \( [1 - \rho, 1 + \rho] \). This yields the perturbed bounds \( \tilde{\mathit{ub}} \), which represent the maximum admissible demand per transport $(i,j)$ and cargo type $k$. Demand samples \( q_{i,j,k} \) are then drawn from a uniform distribution over the interval \( [1, \tilde{\mathit{ub}}_{i,j,k}] \), where the outcome is continuous or discrete depending on whether the model assumes real-valued or integer-valued demand variables. The expected value \( \mu \) and standard deviation \( \sigma \) of the uniform distribution are derived analytically from \( \tilde{\mathit{ub}} \), reflecting the first and second moments of the demand distribution. The resulting tuple \( (q, \mu, \sigma) \) provides necessary information on the demand distribution.

\begin{algorithm}[h]
\caption{Generate Cargo Demand with Uniform Distribution}
\label{alg_drl_scale:uniform_generator}
\textbf{Input:} Capacity $c$, perturbation $\rho$ \\
\textbf{Output:} Realized, expected and standard dev. demand $(q,\mu,\sigma)$ \\
\vspace{0.5em}

Compute upper bound matrix \( \mathit{ub} \propto \mathbf{1}^\top c\)

Apply random perturbation with $ U_{i,j,k} \sim \mathcal{U}(0, 1)$: 
\[
\tilde{\mathit{ub}}_{i,j,k} = \mathit{ub}_{i,j,k}  \times \left(1 + (2 \times U_{i,j} - 1) \times \rho\right) \forall (i, j) \in \textit{TR}, k \in K,
\]

Generate uniform sample: 
\[
q_{i,j,k} \sim 
\begin{cases}
\mathcal{U}(1, \tilde{\mathit{ub}}_{i,j,k}) & \text{if } q \in \mathbb{R}^{n_q}_{>0} \\
\mathcal{U}_{\mathbb{Z}}(1, \tilde{\mathit{ub}}_{i,j,k}) & \text{if } q \in \mathbb{Z}^{n_q}_{>0}
\end{cases}
\quad \forall (i, j) \in \mathit{TR},\; k \in K
\]

Compute expected value: \( \mu = \frac{1}{2} \tilde{\mathit{ub}} \)

Compute standard deviation: \( \sigma = \frac{\tilde{\mathit{ub}}}{\sqrt{12}} \)\\

\Return $q, \mu, \sigma$
\end{algorithm}

\begin{figure}
    \centering
    \begin{subfigure}[b]{0.48\linewidth}
        \centering
        \includegraphics[width=\linewidth]{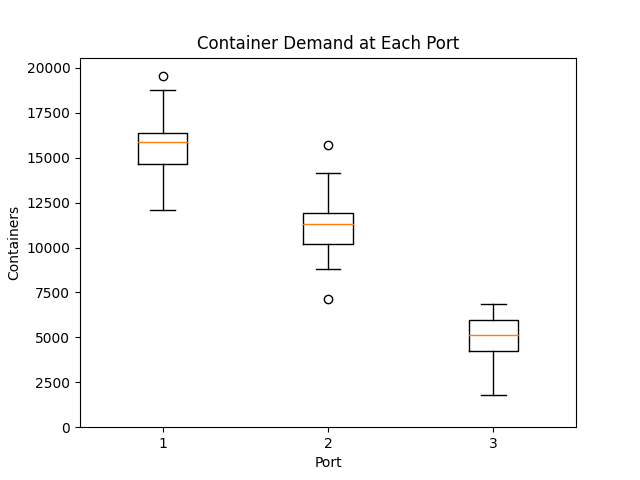}
        \caption{{\small Container demand distribution}}
        \label{fig:containers_p4}
    \end{subfigure}
    \hfill
    \begin{subfigure}[b]{0.48\linewidth}
        \centering
        \includegraphics[width=\linewidth]{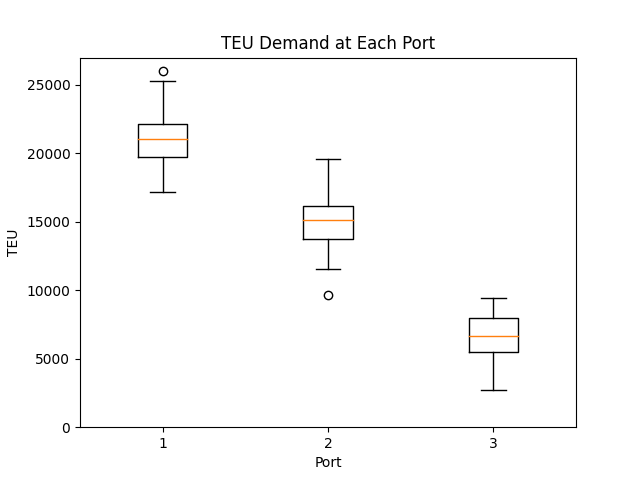}
        \caption{\small TEU demand distribution}
        \label{fig:teu_p4}
    \end{subfigure}
    \caption{Simulated demand for $N_P=4$ by instance generator}
    \label{fig:sim_demand_p4}
    \begin{subfigure}[b]{0.48\linewidth}
        \centering
        \includegraphics[width=\linewidth]{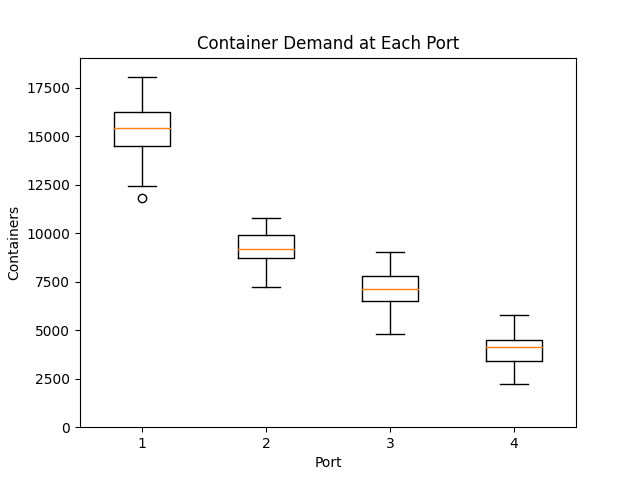}
        \caption{\small Container demand distribution}
        \label{fig:containers_p5}
    \end{subfigure}
    \hfill
    \begin{subfigure}[b]{0.48\linewidth}
        \centering
        \includegraphics[width=\linewidth]{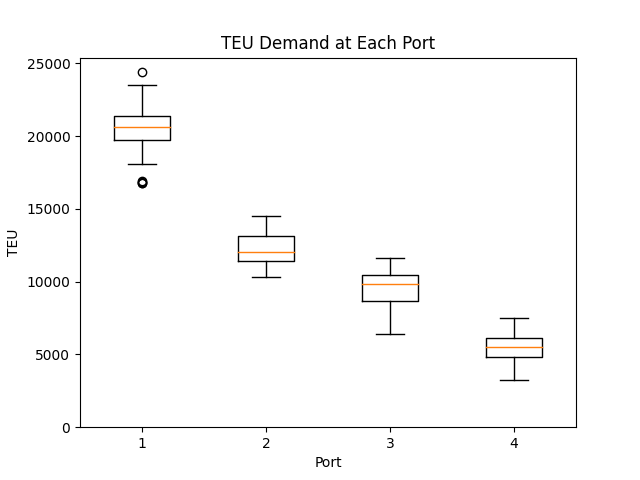}
        \caption{\small TEU demand distribution}
        \label{fig:teu_p5}
    \end{subfigure}
    \caption{Simulated demand for $N_P=5$ by instance generator}
\label{fig:sim_demand_p5}
    \begin{subfigure}[b]{0.48\linewidth}
        \centering
        \includegraphics[width=\linewidth]{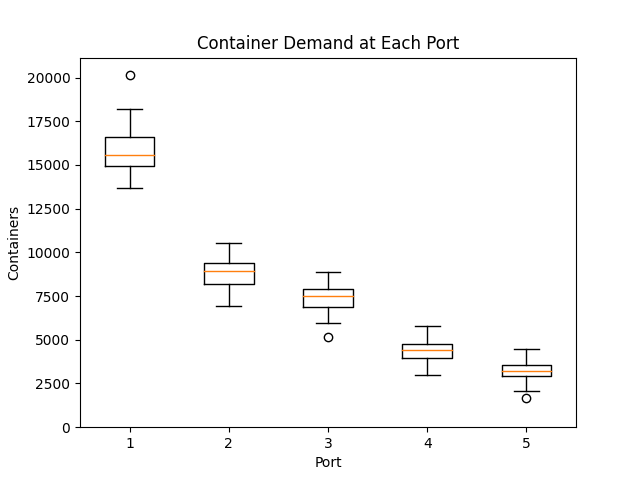}
        \caption{\small Container demand distribution}
        \label{fig:containers_p6}
    \end{subfigure}
    \hfill
    \begin{subfigure}[b]{0.48\linewidth}
        \centering
        \includegraphics[width=\linewidth]{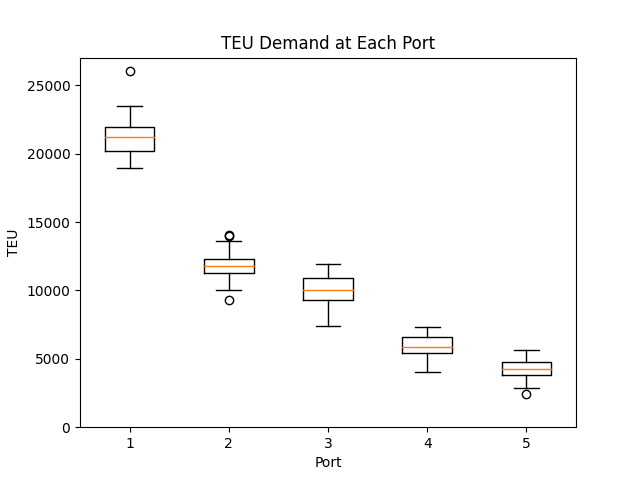}
        \caption{\small TEU demand distribution}
        \label{fig:teu_p6}
    \end{subfigure}
    \caption{Simulated demand for $N_P=6$ by instance generator}
\label{fig:sim_demand_p6}
\end{figure}

\section{Implementation Details of AI2STOW} \label{app_drl_scale:DRL}
In this appendix, we first define the SAC algorithm, then discuss implemented projection layers, and finally describe relevant hyperparameters. Parts of this appendix are adapted from our previous work \cite{van_twiller_navigating_2025}.

\subsection{SAC Algorithm}
Soft Actor-Critic (SAC) is an off-policy reinforcement learning algorithm that optimizes both reward maximization and entropy to encourage efficient exploration \cite{haarnoja_soft_2018}, as outlined in Algorithm \ref{alg_drl_scale:sac}. It is based on maximum entropy reinforcement learning, which aims to learn a stochastic policy that not only maximizes cumulative rewards but also maintains high entropy for robustness and stability. SAC leverages a soft Q-learning approach, using two Q-functions to mitigate overestimation bias, an entropy-regularized policy update, and an automatically adjusted temperature parameter to balance exploration and exploitation.

The algorithm maintains an actor network for policy learning, two Q-function critics for value estimation, a target Q-network for stable learning, and an adaptive temperature parameter to regulate entropy. The loss functions for standard SAC are derived from the Bellman backup equation and the policy gradient formulation, ensuring convergence to an optimal stochastic policy. We also include feasibility regularization from Equation \eqref{for_drl_scale:feas_loss} in the actor loss. 
\begin{itemize}
 \item Compute target Q-value:
\end{itemize}
\begin{align*}
Q_\text{target}(s_t,x_t) &= r_t + \gamma \mathbb{E}_{s_{t+1}, x_{t+1} \sim \pi} \Big[ \\
& \min_{l=1,2} Q_{\theta}^l(s_{t+1}, x_{t+1}) - \alpha \log \pi_{\theta}(x_{t+1} | s_{t+1}) \Big]
\end{align*}
\begin{itemize}
    \item Critic loss:
    \[
    \mathcal{L}_\text{critic}(\theta) = \mathbb{E} \Big[ (Q_{\theta}(s_t, x_t) - Q_\text{target}(s_t,x_t))^2 \Big]
    \]
    \item Actor loss:
    \[
    \mathcal{L}_\text{actor}(\theta) = \mathbb{E} \Big[ \alpha \log \pi_\theta(x_t | s_t) - Q_{\theta}(s_t, x_t) + \lambda_f \mathcal{L}_\text{feas}(\theta)  \Big]
    \]
    \item Temperature loss:
    \[
    \mathcal{L}_\alpha(\theta) = \mathbb{E} \Big[ -\alpha (\log \pi_\theta(x_t | s_t) + \text{entropy}_{\text{target}}) \Big]
    \]

\end{itemize}

This formulation ensures stability and encourages exploration by dynamically adapting the trade-off between exploitation and exploration.

\begin{algorithm}[h!]
\caption{Soft Actor-Critic (SAC)}
\label{alg_drl_scale:sac}
\begin{algorithmic}[1]
\Require Actor parameters $\theta_{\text{actor}}$, critic parameters $\theta_{\text{critic}}^1$, $\theta_{\text{critic}}^2$
\State Initialize target critics: $\theta_{\text{target}}^1 \gets \theta_{\text{critic}}^1$, $\theta_{\text{target}}^2 \gets \theta_{\text{critic}}^2$
\State Initialize temperature parameter $\alpha$
\State Set hyperparameters: learning rates $\eta_a$, $\eta_c$, $\eta_\alpha$; soft update rate $\tau_\textit{su}$; replay buffer $\mathcal{D}$

\For{each training iteration}
    \For{each environment step}
        \State Sample action $x_t \sim \pi_{\theta_{\text{actor}}}(x_t \mid s_t)$
        \State Observe next state $s_{t+1} \sim \mathcal{T}(s_{t+1} \mid s_t, x_t)$
        \State Observe reward $r_t = \mathcal{R}(s_t, x_t)$
        \State Store transition $(s_t, x_t, r_t, s_{t+1})$ in buffer $\mathcal{D}$
    \EndFor

    \For{each gradient update step}
        \State Sample minibatch $(s_t, x_t, r_t, s_{t+1})$ from buffer $\mathcal{D}$

        \State Compute target Q-value:
        \State \quad Sample $x_{t+1} \sim \pi_{\theta_{\text{actor}}}(x \mid s_{t+1})$
        \State \quad $Q_{\text{target}} \gets r_t + \gamma \cdot \min_{l=1,2} Q_{\theta_{\text{target}}^l}(s_{t+1}, x_{t+1}) - \alpha \log \pi_{\theta_{\text{actor}}}(x_{t+1} \mid s_{t+1})$

        \For{$l = 1$ to $2$}
            \State Update critic $\theta_{\text{critic}}^l$ using:
            \State \quad $\theta_{\text{critic}}^l \gets \theta_{\text{critic}}^l - \eta_c \nabla_{\theta_{\text{critic}}^l} \mathcal{L}_{\text{critic}}$
        \EndFor

        \State Update actor parameters:
        \State \quad $\theta_{\text{actor}} \gets \theta_{\text{actor}} - \eta_a \nabla_{\theta_{\text{actor}}} \mathcal{L}_{\text{actor}}$

        \State Update temperature:
        \State \quad $\alpha \gets \alpha - \eta_\alpha \nabla_\alpha \mathcal{L}_\alpha$

        \For{$l = 1$ to $2$}
            \State Soft update target critic:
            \State \quad $\theta_{\text{target}}^l \gets \tau_\textit{su} \cdot \theta_{\text{critic}}^l + (1 - \tau_\textit{su}) \cdot \theta_{\text{target}}^l$
        \EndFor
    \EndFor
\EndFor

\State \Return Policy $\pi_{\theta_{\text{actor}}}$
\end{algorithmic}
\end{algorithm}

\subsection{Projection Layers}
This subsection introduces the projection layers used in the article.

\subsubsection{Violation Projection}
In previous work \cite{van_twiller_navigating_2025}, we discussed the violation projection on a convex polyhedron of constraints. For full technical details, we refer to \cite{van_twiller_navigating_2025}.

Algorithm~\ref{alg:viol_projection} defines a violation projection (VP) layer that reduces inequality constraint violations by shifting a point $x$ closer to the feasible region of the convex polyhedron $\textit{PH} = \{x \in \mathbb{R}^{n}_{>0} : Ax \leq b\}$, where $A \in \mathbb{R}^{m\times n}$ and $b \in \mathbb{R}^{m}$. The inequality constraint $Ax \leq b$ ensures that $\textit{PH}$ is convex, enabling gradient-based reduction of violation $\mathcal{V}(x)$ \cite{boyd_convex_2004}.

To measure feasibility, we compute the element-wise violation $\mathcal{V}(x) = (Ax - b)_{>0}$, where $\mathcal{V}(x)_{m_i} > 0$ indicates that constraint $m_i$ is violated, and $\mathcal{V}(x)_{m_i} = 0$ indicates satisfaction. To minimize violations, we iteratively update $x$ using gradient descent on the squared violation norm $\|\mathcal{V}(x)\|_2^2$, which approximates the distance to the feasible region $\textit{PH}$. The update step is given by:
\begin{align}
x' & = x - \eta_v \nabla_x \|\mathcal{V}(x)\|_2^2   \label{for:x_update_1}
\end{align}
This is equivalent to the update rule shown in Algorithm~\ref{alg:viol_projection}.

During training, the VP layer executes for a fixed number of \texttt{epochs}. However, during inference, a stopping criterion is introduced: the projection halts when the change in total violation, $\mathbf{1}^\top\mathcal{V}(x') - \mathbf{1}^\top\mathcal{V}(x)$, falls below a threshold $\delta_v$. As a result, the VP layer reduces the distance of $x$ to the feasible region, incorporating constraint awareness into otherwise unconstrained policies.

\begin{algorithm}[h!]
\caption{Violation Projection Layer}
\label{alg:viol_projection}
\begin{algorithmic}[1]
\Require Input vector $x \in \mathbb{R}^n_{>0}$, parameters $(A, b, \eta_v, \delta_v)$
\Ensure Projected vector $x'$

\State Initialize $x' \gets x$
\State Define violation function: $\mathcal{V}(x) \gets \max(Ax - b, 0)$ (element-wise)

\For{each iteration $i = 1$ to \texttt{epochs}}
    \State Set $x \gets x'$
    \State Update $x' \gets x - \eta_v A^\top \mathcal{V}(x)$
    \If{$\mathbf{1}^\top \mathcal{V}(x') - \mathbf{1}^\top \mathcal{V}(x) \leq \delta_v$}
        \State \textbf{break}
    \EndIf
\EndFor

\State \Return $x'$
\end{algorithmic}
\end{algorithm}

\subsubsection{Policy Clipping}

In our previous work \cite{van_twiller_navigating_2025}, we also introduced policy clipping (PC) to enforce TEU capacity constraints. Policy clipping applies element-wise lower and upper bounds to a vector $x$, ensuring it remains within a specified box-constrained region. This is implemented using the function:
\[
\mathcal{C}(x,\textit{lb}_\textit{pc},\textit{ub}_\textit{pc}) = \max\big(\min(x, \textit{ub}_\textit{pc}), \textit{lb}_\textit{pc}\big)
\]
where \(\textit{lb}_\textit{pc}\) and \(\textit{ub}_\textit{pc}\) represent the element-wise lower and upper bounds, respectively. While PC is simple and efficient, it is only applicable to box constraints. For more complex convex constraint structures, we rely on the VP layer. Full implementation details are provided in \cite{van_twiller_navigating_2025}.

\subsubsection{Convex Program Layer}
Given a convex polyhedron $\mathit{PH} = \{x \in \mathbb{R}^n \mid Ax \leq b\}$, we incorporate a differentiable convex optimization layer following \cite{agrawal_differentiable_2019}, where the solution to a parameterized convex problem is treated as the output of a network layer. Specifically, we solve
\[
x^\star(\theta) = \arg\min_{x \in \mathit{PH}} f(x, \theta),
\]
where $\theta$ denotes parameters produced by upstream layers and $f(x, \theta)$ is a convex objective. Gradients $\partial x^\star / \partial \theta$ are computed via implicit differentiation of the KKT conditions, enabling end-to-end training. However, since the mapping $\theta \mapsto x^\star(\theta)$ is defined implicitly via the solution of an optimization problem, the Jacobian determinant $\det(\partial x^\star / \partial \theta)$ is not tractable. As a result, the layer cannot be used to compute log-density corrections via the change-of-variables formula, limiting its applicability in likelihood-based or flow-based generative models.
use in likelihood-based models.

To promote numerical stability, we incorporate slack variables into the stability constraint and apply a large penalty (e.g., $10^4$) for violations. Specifically, we relax constraints of the form $g(x) \leq 0$ to $g(x) \leq \varepsilon$, and penalize $\varepsilon$ in the objective via $\lambda_\textit{cp} \|\varepsilon\|_1$. This soft enforcement accommodates intermediate actions that may temporarily breach stability constraints, as can arise in the MDP.

\newpage
\subsection{Hyperparameters}
\begin{table*}[h!]
    \centering
    \small
    \caption[Hyperparameters AI2STOW]{Hyperparameters of AI2STOW. Note that the hyperparameters are identical for DRL-FR, except for $\lambda_f$ as we performed parameter tuning to set it to $\lambda_f=0.195$} \label{tab:ppo_sac_hyperparameters}
    \begin{tabular}{llc}
        \toprule
        \textbf{Hyperparameters} & \textbf{Symbol} & \textbf{AI2STOW} \\
        \midrule
        \textbf{Actor Network} & & Attention \\
        \textbf{Number of Heads} & & 8  \\
        \textbf{Hidden Layer Size} & & 512  \\
        \textbf{Encoder Layers} & & 3 \\
        \textbf{Decoder Layers} & & 4 \\
        \textbf{Critic Network} & & $2\times \text{MLP}$ \\
        \textbf{Critic Layers} & & 5\\
        \textbf{Target Network} & & Soft Update \\
        \textbf{Target Update Rate} & $\tau_\textit{su}$ & 0.005 \\
        \textbf{Dropout Rate} & & $0.009$\\
        \textbf{Max Policy Std.} & & $9.460$ \\
        \midrule
        \textbf{Optimizer} & & Adam\\
        \textbf{Learning Rate} & $\eta$ & $1.46 \times 10^{-4}$ \\
        \textbf{Batch Size} & & 64 \\
        \textbf{Embedding Size} & & 128 \\
        \textbf{Discount Factor} & $\gamma$ & 0.99 \\
        \textbf{Entropy Coefficient} & $\lambda_e$ & Learned \\
        \textbf{Feasibility Penalty} & $\lambda_f$ & $0.283$ \\
        \textbf{Replay Buffer Size} & &  $10^4$ \\
        \textbf{Mini-batch Size} & & 32 \\
        \textbf{Entropy Target} & & $-|{X}|$\\
        \midrule
        \textbf{Projection Learning Rate} & $\eta_v$ & $0.010$ \\
        \textbf{Projection Epochs} & & $273$ \\
        \textbf{Inference Projection Stop} & $\delta_v$ & $0.024$\\
        \textbf{Finetuned Projection Learning Rate} & $\eta_v$\textsuperscript{\ding{72}} & $0.01$ \\
        \textbf{Finetuned Projection Epochs} & & $300$ \\
        \textbf{Finetuned Inference Projection Stop} & $\delta_v$\textsuperscript{\ding{72}} & $0.01$\\
        \textbf{Slack penalty of CP} & $\lambda_\textit{cp}$ & $1 \times 10^4$\\
        \midrule
        \textbf{Training Budget} & & $7.2 \times 10^{6}$\\
        \textbf{Validation Budget} & & $5.0 \times 10^{3}$\\
        \textbf{Validation Frequency} & & Every 20\% \\
        \textbf{Inference Rollouts} & & $5$\\
        \bottomrule
    \end{tabular}
\end{table*}



\section{Additional Experiments} \label{app_drl_scale:add_experiments}

\Cref{tab_drl_scale:memory} illustrates the rapid growth in computational demands of the scenario tree in the SMIP-NA model with $S_\textit{ST} = 20$. As the voyage length $N_P$ increases, the number of scenario paths $|\mathcal{Z}|$, as well as the runtime and memory usage, grow substantially. For $N_P = 4$, both runtime and memory consumption are measured directly and serve as a baseline. Assuming an optimistic case of linear scaling with respect to $|\mathcal{Z}|$, the extrapolated values for $N_P = 5$ and $N_P = 6$ already exceed practical hardware and runtime limits by a significant margin. This highlights the challenge of solving large-scale instances of SMIP-NA with conventional computational resources.

\begin{table*}[h!]
\centering
\small
\caption[Evaluation of memory use scenario tree stochastic MIP model]{Evaluation of runtime and peak memory usage for the SMIP-NA model with $S_\textit{ST} = 20$, as the scenario tree expands with increasing voyage length $N_P$. Runtime and memory usage for $N_P = 5$ and $N_P = 6$ are extrapolated linearly based on the growth of $|\mathcal{Z}|$ relative to $N_P = 4$, and are marked with an asterisk ($\ast$) to indicate estimated values.}
\label{tab_drl_scale:memory}
\begin{tabular}{lrrr}
\toprule
\multicolumn{1}{c}{\textbf{SMIP-NA}} & \multicolumn{3}{c}{\textbf{Metrics}} \\
\cmidrule(r){1-1} \cmidrule(r){2-4}
\textbf{$N_P$} & \textbf{ $|\mathcal{Z}|$} & \textbf{Runtime (s)} & \textbf{Memory (GB)} \\
\midrule
4 & 400 & 1576.77 &  36.01 \\
5 & 8,000 &  31,535.40\textsuperscript{*} &  720.23\textsuperscript{*} \\
6 & 160,000 &  630,708.00\textsuperscript{*} &  14,404.66\textsuperscript{*} \\
\bottomrule
\end{tabular}
\end{table*}


\end{document}